\documentclass[12pt]{article}
\usepackage{amssymb,amsmath,makeidx,amsthm}
\usepackage{enumerate}
\usepackage{psfig}
\usepackage{color}
\usepackage[all]{xy}
\begin{document}
\newtheorem{theorem}{Theorem}[subsection]
\newtheorem{pusto}[theorem]{}
\numberwithin{subsection}{section}
%\numberwithin{figure}{chapter}
\newtheorem{corollary}[theorem]{Corollary}
\newtheorem{remark}[theorem]{Remark}
%\newtheorem{example}[theorem]{Example}
%\newtheorem{exercise}[theorem]{Exercise}
%\newtheorem{conjecture}[theorem]{Conjecture}
%\newtheorem{question}[theorem]{Question}
%\newtheorem{problem}{Problem}
%\newtheorem{criterion}[theorem]{Criterion}
%\newtheorem{generalization}[theorem]{Generalization}
%\newtheorem{conjecture}[problem]{Conjecture}
%\input newmacro.tex
%\input mssymb12.tex
%%%%%%%%%%%%%%%%%%%%%%%%%%%%%%%%%%%%%%%%%%%%%%%%%%%%%%%%%%%
%\newcommand{\proof}[1]{\noindent{\bf Proof#1:\  }}
%\newcommand{\qed}{\hfill$\square$\medskip}
\newcommand{\po}{{\perp_{\omega}}}
\newcommand{\mathto}{\mathop{\longrightarrow}\limits}
%%%%%%%%%%%%%%%%%%%%%%%%%%%%%%%%%%%%%%%%%%%%%%%%%%%%%%%%%%%
%\def\bt{\begin{theorem}}
%\def\et{\end{theorem}}
\def\x{\index}
\def\bp{\begin{pusto}}
\def\ep{\end{pusto}}
\def\bi{\begin{itemize}}
\def\ei{\end{itemize}}
\def\btl{$\blacktriangleleft$}
\def\btr{$\blacktriangleright$}
\def\mn{\medskip\noindent}
\def\n{\noindent}
\def\p{\partial}
\def\e{\varepsilon}
\def\a{\alpha}
\def\d{\delta}
\def\s{\sigma}
\def\ol{\overline}
\def\wt{\widetilde}
\def\Cal{\cal}
%%%%%%%%%%%%%%%%%%%%%%%%%%%%%%%%%%%%%%
\def\Int{\mathrm{Int}}
\def\Ker{\mathrm{Ker}}
\def\Id{\mathrm{Id}}
\def\id{\mathrm{id}}
\def\su{\mathrm{supp}}
\def\rank{\mathrm{rank}}
\def\dim{\mathrm{dim}}
\def\codim{\mathrm{codim}}
\def\const{\mathrm{const}}
\def\Span{\mathrm{Span}}
\def\span{\mathrm{span}}
\def\inv{\mathrm{inv}}
\def\Hom{\mathrm{Hom}}
\def\Mono{\mathrm{Mono}}
\def\Diff{\mathrm{Diff}}
\def\Emb{\mathrm{Emb}}
\def\emb{\mathrm{emb}}
\def\imm{\mathrm{imm}}
\def\sub{\mathrm{sub}}
\def\Hol{\mathrm{Hol}}
\def\hol{\mathrm{hol}}
\def\Sec{\mathrm{Sec}}
\def\Sol{\mathrm{Sol}}
\def\Supp{\mathrm{Supp}}
\def\bs{\mathrm{bs}}
\def\Vect{\mathrm{Vect}}
\def\Vert{\mathrm{Vert}}
\def\Clo{\mathrm{Clo}}
\def\clo{\mathrm{clo}}
\def\Exa{\mathrm{Exa}}
\def\Gr{\mathrm{Gr}}
\def\CR{\mathrm{CR}}
\def\Distr{\mathrm{Distr}}
\def\dist{\mathrm{dist}}
\def\Ham{\mathrm{Ham}}
\def\mers{\mathrm{mers}}
\def\Iso{\mathrm{Iso}}
\def\iso{\mathrm{iso}}
\def\tang{\mathrm{tang}}
\def\trans{\mathrm{trans}}
\def\ttr{\mathrm{trans-tang}}
\def\ot{\mathrm{ot}}
\def\loc{\mathrm{loc}}
\def\Symp{\mathrm{Symp}}
\def\symp{\mathrm{symp}}
\def\Lag{\mathrm{Lag}}
\def\isosymp{\mathrm{isosymp}}
\def\subisotr{\mathrm{sub-isotr}}
\def\Cont{\mathrm{Cont}}
\def\cont{\mathrm{cont}}
\def\Leg{\mathrm{Leg}}
\def\isocont{\mathrm{isocont}}
\def\isot{\mathrm{isot}}
\def\isotr{\mathrm{isotr}}
\def\coisot{\mathrm{coisot}}
\def\coreal{\mathrm{coreal}}
\def\real{\mathrm{real}}
\def\comp{\mathrm{comp}}
\def\old{\mathrm{old}}
\def\new{\mathrm{new}}
\def\top{\mathrm{top}}
\def\rG{\mathrm{G}}
\def\CS{\mathrm{CS}}
\def\ND{\mathrm{ND}}
\def\CND{\mathrm{CND}}
%%%%%%%%%%%%%%%%%%%%%%%%%%%%%%%%%%%%%%%%%%%%%%
\def\fA{\mathfrak {A}}
\def\fa{\mathfrak {a}}
\def\A{\mathcal {A}}
\def\E{\mathcal {E}}
\def\F{\mathcal {F}}
\def\G{\mathcal {G}}
\def\H{\mathcal {H}}
\def\I{\bf {I}}
\def\J{\mathcal {J}}
\def\N{\mathcal {N}}
\def\R{\mathcal {R}}
\def\P{\mathcal {P}}
\def\T{\mathcal {T}}
\def\L{\mathcal {L}}
\def\SS{\mathcal {S}}
\def\U{\mathcal {U}}
\def\D{\mathcal {D}}
\def\V{\mathcal {V}}
\def\W{\mathcal {W}}
\def\Z{\mathcal {Z}}
\def\C{\mathcal {C}}
\def\O{\mathcal {O}}
%%%%%%%%%%%%%%%%%%%%%%%%%%%%%%%%%%%%%%%%%%%%%%%%%
\def\BH{{\bold H}}
\def\BF{{\bold F}}
\def\BS{{\bold S}}
\def\bbZ{{\mathbb Z}}
\def\bbR{{\mathbb R}}
\def\bbS{{\mathbb S}}
\def\bbC{{\mathbb C}}
\def\bbQ{{\mathbb Q}}
%%%%%%%%%%%%%%%%%%%%%%%%%%%%%%%%%%%%%%%%%%%%%%%%%
\def\IS{\R_{\sympl-iso}}
\def\IC{\R_{\cont-iso}}
\def\Sp{\mathrm{Sp}}
\def\GL{\mathrm{GL}}
\def\SO{\mathrm{SO}}
\def\U{\mathrm{U}}
\def\X{X^{(r)}}
\newcommand{\eps}{\epsilon}
\def\Op{{\mathcal O}{\it p}\,}
\def\hook{\lrcorner\,}

\newcommand{\corank}{\mathrm{Corank}}
\title{
WRINKLED EMBEDDINGS \\
{\small To Paul Schweitzer on his 70th birthday}}
\author{
     Y. M. Eliashberg
     \thanks{Partially supported by the NSF grant DMS-0707103}
     \\Stanford University, \\Stanford, CA 94305 USA
\and
 N. M. Mishachev
%      \thanks{Partially supported by the NSF }
  \\ Lipetsk Technical University, \\
     Lipetsk, 398055 Russia}

\date{}
\maketitle

\begin{abstract}
\noindent
A {\it wrinkled embedding} $f:V^n\to W^m$
is a topological embedding which is a smooth embedding
everywhere on $V$ except a set
of $(n-1)$-dimensional spheres, where $f$ has cuspidal corners.
In this paper we prove that any rotation of the tangent
plane field $TV\subset TW$ of a {\it smoothly embedded} submanifold
$V\subset W$
can be approximated by a homotopy of {\it wrinkled embeddings} $V\to W$.
\end{abstract}

\tableofcontents
%\clearpage
\setcounter{section}{0}

\section{Introduction}\label{intro}

\subsection{Wrinkled embeddings}\label{sec:wr-embed-informal}
\begin{figure}[hi]
\centerline {\psfig{figure=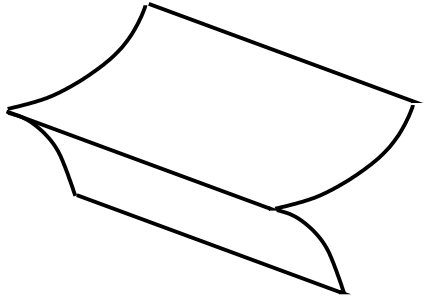,height=40mm}}
\caption {\small The wrinkled embedding $f:V^2\to W^3$
near $S_i\setminus S'_i$}
\label{we01}
\end{figure}

\begin{figure}[hi]
\centerline{\psfig{figure=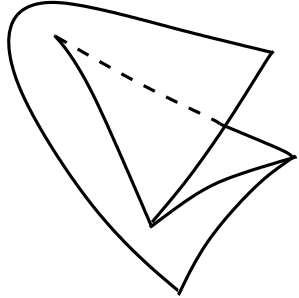,height=40mm}}
\caption {\small The wrinkled embedding $f:V^2\to W^3$
near $S'_i$}
\label{we02}
\end{figure}

A {\it wrinkled embedding} $f:V^n\to W^m$, $n<m$,
is a topological embedding which is a smooth embedding
everywhere on $V$ except a finite set
of $(n-1)$-dimensional spheres $S_i\subset V$, where $f$ has cuspidal corners:
threefold corners along an equator $S'_i\subset S_i$ and twofold corners
along the complement $S_i\setminus S'_i$. For $n=2$ and $q=3$ see
Fig.\,\ref{we01} and Fig.\,\ref{we02}. The spheres $S_i$ are called 
{\it wrinkles}. A formal definition of wrinkled embeddings is given in Section
\ref{sec:wr-embed-formal} below.

\n
Note that for $n=1$ each wrinkle consists of two
twofold cuspidal points.
{\it Families} of wrinkled embeddings may have,
in addition to wrinkles, {\it embryos} of
wrinkles, and therefore wrinkles may {\it appear} and {\it disappear}
in a homotopy of  wrinkled embeddings.

Let us point out that  a {\it wrinkled embedding is not a wrinkled map}
in the sense of \cite{[EM97]}.
The relation between the two notions is discussed in Section
\ref{sec:wr-embed-formal} below.

\subsection{Main result and  the idea of the proof}\label{sec:main}

\begin{figure}[hi]
\centerline{\psfig{figure=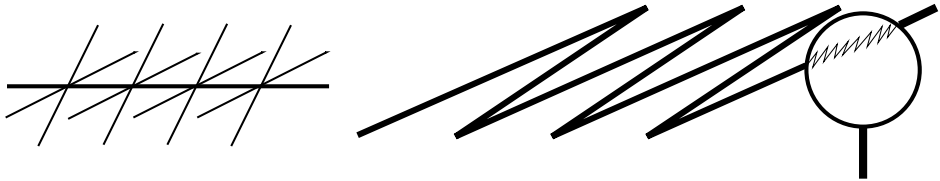,height=25mm}}
\caption {\small Wrinkled approximation of a counterclockwise tangential
rotation}
\label{we03}
\end{figure}

In this paper we prove (see Theorem  \ref{thm:int-tang-global} below)
that any homotopy of the tangent
plane field $TV\subset TW$ ({\it tangential homotopy})
of a {\it smoothly embedded} submanifold $V\subset W$
can be approximated by a {\it wrinkled  isotopy}, i.e. an isotopy through
wrinkled embeddings $V\to W$.
For $n=1$ the idea of the proof is presented
on  Fig.\,\ref{we03}. Here we consider the counterclockwise tangential
rotation of an interval in $\bbR^2$.

\mn
The implementation of this general idea for $n>1$ is far from being
straightforward. There are lot of similarities here with  the Nash-Kuiper
theorem about isometric $C^1$-embeddings $V^n\to W^{n+1}$,
where the proof in the case $n=1$ is more or less trivial,
while already contains the general idea (goffering).
However, its realization for $n>1$ is highly non-trivial.

\subsection{Applications}\label{sec:applic}

Our main theorem can be reformulated as an
{\it $h$-principle for $A$-directed wrinkled embeddings}, see
\ref{thm:dir-wrinkled}.
As an application of the main theorem  we prove an $h$-principle
for {\it embeddings\,} $f:V^n\to (W,\cal \xi)$, $n\geqslant q=\codim\,\xi$,
whose tangency  singularities with respect to a
distribution $\xi$ (integrable or non-integrable)
are simple, e.g. folds, or alternatively  {\it generalized} wrinkles
(see Theorem \ref{thm:embedding-foliation} and Section \ref{sec:embed-to-dis}).
This $h$-principle allows us, in particular, to simplify
the singularities of an individual embedding or a family of embeddings
$V^n\to \bbR^m$, $n<m$,
with respect to the projection
$\bbR^m=\bbR^q\times \bbR^{m-q}\to\bbR^q$, $n\geqslant q$,
i.e with respect to the standard foliation of $\bbR^m$ by $(m-q)$-dimensional
affine subspaces, parallel to $0\times\bbR^{m-q}$.

\subsection{History of the problem}\label{sec:hist}

{\bf A. Directed embeddings.}
An embedding $V^n\to W^m$ is called {\it $A$-directed}, if its tangential
(or Gaussian) image belongs to a given subset $A$ of the Grassmannian
bundle $\Gr_nW$.
Using his convex integration method, Gromov proved in \cite{[Gr86]}
a general $h$-principle for $A$-directed embedding in the case when $V$ is
an {\it open} manifold and $A\subset \Gr_nW$ is an open subset.
C.~Rourke and B.~Sanderson gave two independent proofs of this
theorem in \cite{[RS01]}.  A different proof based on our
holonomic approximation theorem was given in our book \cite{[EM02]}.
For some special $A\subset\Gr_nW$ Gromov  also proved in \cite{[Gr73]} and
\cite{[Gr86]}  the $h$-principle for {\it closed} manifolds.
However, for {closed} $V$ the $h$-principle for $A$-directed embeddings fails
for a {\it general} open $A$. For example, for any closed $V^n$ there is
no $A$-directed embeddings $V^{n}\to \bbR^{n+1}$ unless $\pi(A)=S^n$, where
$$\pi:\Gr_n\bbR^{n+1}=S^n\times\bbR^{n+1}\to S^n$$
is the projection.
The main theorem  of the current paper states that
{\it this $h$-principle can be saved by relaxing the notion of embedding}.
D.~Spring in \cite{[Sp05]} proved, using Gromov's convex integration method
and the {\it geometry of spiral curves} from \cite{[Sp02]},
an existence theorem for directed embeddings with {\it twofold spherical corners},
which is equivalent to the non-parametric version of our main theorem,
see the discussion in Sections  \ref{sec:double} and \ref{sec:direct} below.

\mn
{\bf B. Simplification of singularities.}
%This paper continues our {\it wrinklization} program,
%see  \cite{[EM97]},\cite{[EM98]},\cite{[EM00]}.
First result allowing to simplify singularities of an {\it  individual}
embedding $V^n\to \bbR^m$, $n<m$, with respect to the projection
$\bbR^m=\bbR^q\times \bbR^{m-q}\to\bbR^q$, $n\geqslant q$,
was proven in \cite{[El72]}. In fact, in \cite{[El72]} this result was formulated
for $n=m-1$; the general case $n<m$ can be derived
from this basic one by  induction which was done
  in \cite{[EM00]}. A different proof based on
convex integration method was given by D.~Spring in \cite{[Sp02]}.
In \cite{[EM00]} we also proved the {\it parametric} version,
but only for $q=1$, and only the epimorphism
part of the corresponding parametric $h$-principle.
However,  the approach in \cite{[EM00]}  does not seem to be suitable
to recover the main results of the current paper.
%\mn
%{\bf C. Problems.}
%Let us remind, that the Smale conjecture $\Diff(S^3)\simeq {\rm O}(4)$
%is equivalent to the

The reader   may find   additional interesting information
related to the subject of this paper  in \cite{[Sp02]} and \cite{[Sp05]}.
A different  approach to  the problem of simplification of singularities 
can be found in \cite{[RS03]}.

\subsection{Remark}\label{sec:rem}
We assume that the reader is familiar  with the general philosophy
of the $h$-principle, see \cite{[Gr86]} and \cite{[EM02]}.
It is useful, though not necessary for  the reader to be  also  familiar with
the Holonomic Approximation Theorem  from \cite{[EM02]}) and  the wrinkling
philosophy (see  \cite{[EM97]} and \cite{[EM98]}).
%Though formally the wrinkling {\it theorem} does not play any role in our proofs,
We recall in Section \ref{appendix}, for a convenience of the reader,
some definitions and results from \cite{[EM97]} and \cite{[EM98]}
and introduce there the notions of {\it generalized wrinkles} and
{\it generalized wrinkled maps}.

\section[Integrable approximations of tangential
homotopies] {Integrable approximations of tangential\\
homotopies}\label{sec:integral-approx}

\subsection{Tangential homotopies of embeddings}\label{sec:tang-homot}

In what follows we assume that
$V\subset W$ is an embedded compact submanifold and denote by
$f_0$ the inclusion ${\rm{i}}_V:V\hookrightarrow W$.
We also assume that the manifolds $W$ and $\Gr_nW$ are endowed with
Riemannian metrics.

\mn
Let $\pi:\Gr_nW\to W$
be the Grassmannian bundle of  tangent $n$-planes
to a $m$-dimensional manifold $W$, $m>n$, and $V$ a $n$-dimensional
manifold. Given a {\it monomorphism}
(= fiberwise injective homomorphism) $F:TV\to TW$,
we will denote by $\rG F$ the corresponding map
$V\to \Gr_nW$. Thus the tangential (Gaussian) map associated
with an immersion $f:V\to W$ can be written as $\rG df$.

\mn
A {\it tangential homotopy} of
an embedding $f_0$ is a  homotopy $G_t:V\to \Gr_nW$, such that $\,G_0=\rG df_0$
which covers an isotopy  $g_t=\pi\circ G_t:V\to W$, $g_0=f_0$.
A tangential homotopy  $G_t:V\to \Gr_nW$
is called {\it integrable}, if $G_t=\rG dg_t$.
A tangential homotopy $G_t:V\to \Gr_nW$ is called a
{\it tangential rotation}, if $\,G_0=\rG df_0$ and $\pi\circ G_t=f_0$,
i.e. $G_t$ covers the constant isotopy $g_t=f_0:V\to W$.

\mn
\btl\,
{\bf Problem.}
Let $G_t:V\to \Gr_nW$ be a tangential homotopy.
We want to construct an arbitrarily close
{\it integrable approximation} of $G_t$, i.e
an isotopy of embeddings $f_t:V\to W$, such that  $\rG df_t$ is arbitrarily
close to $G_t$.
One can reduce the problem to the case
when $G_t$ is a {\it tangential rotation}. Indeed, we can
consider, instead of $G_t$, the rotation $(d\hat g_t)^{-1}\circ\, G_t$,
where $\hat g_t: W\to W$ is a diffeotopy which extends the
isotopy $g_t=\pi\circ G_t:V\to W$.
Note that if $G_t$ is an {\it integrable}
tangential homotopy, then $(d\hat g_t)^{-1}\circ\, G_t$ is
the constant homotopy $G_t=G_0$.
\btr

\mn
Of course,   an integrable approximation of a tangential rotation
does not exist in general.
However, as we shell see, one can always construct an
integrable approximation of $G_t$ by a family of
{\it wrinkled embeddings}.

\subsection{Local integrable approximations}\label{sec:local-integr}

Let $X\subset W$ be a tubular neighborhood of $V\subset W$,
and $\pi:X\to V$ the normal projection. Denote by $\N$ the
normal foliation of $X$ by the fibers of $\pi$.
An isotopy $f_t:V\to W$, $f_0=i_V$,
is called {\it graphical}, if all the images $f_t(V)$ are transversal to $\N$.
In other words, the graphical isotopy $V\to W$
is an isotopy of sections $V\to X$, up to reparameterizatons of $V$.
A tangential rotation $G_t:V\to \Gr_nW$ is called {\it small},
if $G_t(v)$ is transversal to $\N$ for all $t$ and $v$.

\mn
Following Gromov's book \cite{[Gr86]} we will use the
notation $\Op A$ as a replacement of the expression
{\it an open neighborhood of $A\subset V$}. In other words, $\Op A$
is an arbitrarily small but not specified open neighborhood
of a subset $A\subset V$.
\bp
\label{thm:int-tang-local}
{\bf (Local integrable approximation of  a small tangential rotation)}
Let $G_t:V\to \Gr_nW$ be a small tangential rotation of $V\subset W$
and $K\subset V$ a stratified subset of {\bf positive} codimension.
Then there exists an arbitrarily $C^0$-small
graphical isotopy $f_t:V\to W$
such that the homotopy
$$\rG df_t|_{\Op K}:\Op K\to \Gr_nW$$
is arbitrarily $C^0$-close to the tangential rotation $G_t|_{\Op K}$.
\ep

\n
{\bf Proof.}
The space $X^{(1)}$ of
$1$-jets of sections $V\to X$ can be interpreted as a space of
tangent to $X$ $n$-planes which are non-vertical, i.e. transversal
to $\N$. Hence the inclusion
$f_0:V\hookrightarrow X$ together with the tangential homotopy
$G_t:V\to \Gr_nW$ can be viewed  as a homotopy of sections
$F_t:V\to X^{(1)}$.
For arbitrarily small $\e$ and $\d$ we can construct, using
Holonomic Approximation  Theorem
3.1.2 in \cite{[EM02]}, a holonomic $\e$-approximation
$\wt F_t$ of $F_{t}$ over $\Op \wt h_t(K)$,
where $\wt h_t:V\to V$ is a $\d$-small
diffeotopy. The $0$-jet part $\wt f_t$ of the
section $\wt F_t$ is automatically an embedding, because $\wt f_t$ is a section
of the normal bundle. Thus, we have a family $\wt f_t$ of
integrable approximations of $G_t$ over $\Op \wt h_t(K)$,
and hence one can define the required isotopy $f_t:V\to W$
as the composition $\wt f_t\circ \wt h_t$,
where $\wt f_t:V\to X$ is an extension
of the isotopy $\wt f_t:\Op \wt h^t(K)\to X$ to an
isotopy of {\it sections} $V\to X$.
\qed

\n
\btl\,
{\bf Remarks.\,}
\n
{\bf 1.}
One can also apply the above  construction  to {\it any} section $\wt V\subset X$
instead of  the zero-section $V\subset X$, provided that
the tangential rotation of $T\wt V$ is transversal to $\N$
(here $\N$ is the {\it original} foliation on $X$ by the fibers of
the projection $\pi: X\to V$).

\n
{\bf 2.} Without the ``graphical" restriction on $f_t$
Theorem \ref{thm:int-tang-local} remains true for {\it any}
tangential rotation, see 4.4.1 in \cite{[EM02]}.

\n
{\bf 3.} The relative and the parametric versions  of the holonomic
approximation theorem similarly prove
the relative and the parametric versions of
Theorem \ref{thm:int-tang-local}. In the relative
version the homotopy $G_t$ and the diffeotopy $h_t$ are constant near
a compact subset $L\subset K$. In the parametric version
we deal with a family of embedded submanifolds
$V_p\hookrightarrow W$, $p\in B$.  The
parameters space $B$ usually is a manifold, possibly  with non-empty
boundary $\p B$, and the homotopies $G_{t,p}$\,, $p\in B$,
are constant for $p\in \Op \partial B$. In this case
the diffeotopies $h_{t, p}\,$,
$p\in B$, are constant for $p\in \Op \partial B$.
\btr

\subsection{Wrinkled embeddings}\label{sec:wr-embed-formal}

{\bf A. Definition.}
A (smooth) map $f:V^n\to W^m\,,\,\,n<m\,,$
is called  a {\it wrinkled embedding}, if
\bi
\item $f$ is a topological embedding;
\item any connected component $S_i$ of the singularity
$S=\Sigma(f)$ is diffeomorphic
to the standard $(n-1)$-dimensional sphere $S^{n-1}$
and bounds an $n$-dimensional disk $D_i\subset V$;
\item the map $f$ near each sphere $S_i$ is
equivalent to the map
$$Z(n,m):\Op_{\bbR^n} S^{n-1}\to \bbR^{m}$$
given by the formula
 $$(y_1,...,y_{n-1},z)\longmapsto$$
$$ \left(y_1,...,y_{n-1},z^3+3(|y|^2-1)z,
\int_0^z (z^2+|y|^2-1)^2dz,\,0,...,0\right).\;\; \eqno(1)
$$
\ei

\mn
The compositions of the canonical form (1) with the projection
to the space of first $n$ local coordinates in $\bbR^m$ is
the standard equidimensional wrinkled map
$\bbR^n\to \bbR^n$, see Section \ref{sec:wrinkles} or \cite {[EM97]}.
Thus, the canonical form
for the wrinkled embeddings  contains, in comparison to the canonical form for
the wrinkled mappings, the {\it unfolding} function
$\int_0^z (z^2+|y|^2-1)^2\,dz$
and zero functions as additional coordinates in the
image.
The spheres $S_i$ and its images $f(S_i)$ are called {\it wrinkles} of the
wrinkled embedding $f$.
According to the formula (1), each wrinkle  $S_i$
has a marked $(n-2)$-dimensional equator $S'_i\subset S_i$ such that
\bi
\item the local model for $f$ near each point of $S_i\setminus S'_i$ is
given by the formula
$$(y_1,...,y_{n-1},z)\longmapsto
(y_1,...,y_{n-1},z^2, z^3, 0,...,0)\,,\eqno (2)$$
see Fig.\,\ref{we01};
\item the local model for $f$ near each point of $S'_i$ is
given by the formula
$$(y_1,...,y_{n-1},z)\longmapsto
(y_1,...,y_{n-1},z^3-3y_{1}z, \int_0^z (z^2-y_{1})^2dz,0,...,0),\eqno (3)$$
see Fig.\,\ref{we02}.
\ei

\n
The compositions of the canonical forms (2) and (3) with the projection
to the space of first $n$ local coordinates in $\bbR^m$ are
the standard fold and cusp of the equidimensional map
$\bbR^n\to \bbR^n$.

We will denote the union  $\bigcup\limits_i S_i'$ of
all the equators $S'_i$ by $S'$.

\mn
The restriction of  a wrinkled embedding
$f$ to $\Int\,D_i$ {\it may have singularities},\,
which are, of course,
again wrinkles. We will say, that a wrinkle $S_i$ {\it has a
depth} $d+1$, if $S_i$ is contained in $d$ ``exterior" wrinkles.
The depth of a wrinkled embedding is the maximal depth of its wrinkles;
according to this definition the smooth embeddings are the
wrinkled embeddings of depth $0$.

\btl\,
{\bf Remark.}
Note that  a wrinkled embedding  
$f$ near $S\setminus S'$ is equivalent to the restriction
of a generic $\Sigma^{11}$-map of an $(n+1)$-dimensional manifold to its fold
$\Sigma^1$  near the cusp $\Sigma^{11}$.
However, near $S'$  a wrinkled embedding $f$ {\it is not equivalent}
to the restriction of a generic $\Sigma^{111}$-map
to $\Sigma^1$ near the swallow tail singularity $\Sigma^{111}$. \btr

\begin{figure}[hi]
\centerline{\psfig{figure=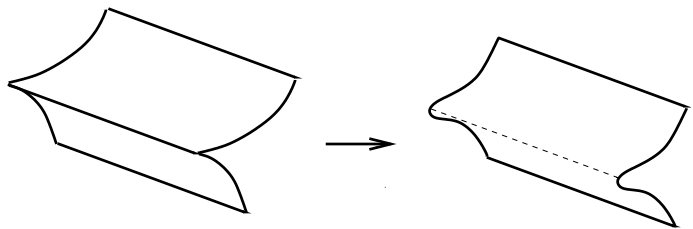,height=30mm}}
\caption {\small Regularization of a wrinkled
embedding near $\Sigma_i\setminus\Sigma'_i$}
\label{we06}
\end{figure}

\begin{figure}[hi]
\centerline{\psfig{figure=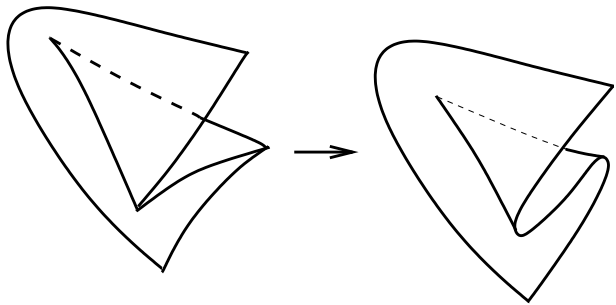,height=40mm}}
\caption {\small Regularization of a wrinkled embedding near $\Sigma'_i$}
\label{we07}
\end{figure}

\mn
{\bf B. Regularization.}
Any wrinkled embedding can be canonically {\it regularized}
by changing the unfolding function $u(y,z)=\int (z^2+|y|^2-1)^2\,dz$
in the canonical form to a $C^1$-close function $\wt u(y,z)$ such that
$\partial_z\wt u(y,z)>0$, see Fig.\,\ref{we06} and Fig.\,\ref{we07}.
One can chose $\wt u$ such that $\wt u(y,z)=u(y,z)$ for all 
$(y,z)\in S\setminus \Op S'$. Then the respective regularization does not moves
the twofold corner points $f(S\setminus S')$ everywhere except an arbitrarily 
small neighborhood of the threefold corner points  $f(S')$.  

%the (i.e. the function $\int (z^2+|y|^2-1)^2\,dz$)
%to the function $z$,
%This is true also in the parametric case (see {\bf D}).

\mn
{\bf C. Gaussian map.}

For any wrinkled embedding $f:V\to W$ and any $v_0\in \Sigma_f$
there exists a limit
$$\lim_{v\to v_{_0}, \,\,v\in V\setminus\Sigma_f}\, f_*(T_vV).$$
Hence, we can associate with a wrinkled embedding its
{\it wrinkled tangent bundle} $T(f)$. If $V$ is oriented then the bundle
$T(f)$ is   oriented as well. The orientation
of a tangent space $T_v$ to $V$  at a smooth point $v\in V\setminus\Sigma(f)$ as
a fiber of $T(f)$  coincides with,  or opposite  to its orientation as the fiber
of $TV$,  depending on the parity of
of the number of exterior wrinkles which surround $v$.
Therefore one  can canonically define a Gaussian map
$\rG df:V\to \Gr_nW$ into the Grassmannian of  non-oriented $n$-planes,
and  if $V$ is oriented   the map  $\wt\rG df:V\to\widetilde{\Gr}_nW$ into
the Grassmannian $\widetilde{\Gr}_nW$ of {\it oriented}
$n$-planes in $TW$.

\mn
{\bf D. Normal foliation.}
Though the image the  $f(V)\subset W$ of a wrinkled embedding
$f:V\to W$  is not a smooth submanifold, one can still define
 an analog of normal foliation.

\mn
There exists an $n$-dimensional submanifold
$\widehat V\subset W$, such that $f(\Sigma)\subset \widehat V$ and
$\widehat V$ is tangent to $f(V)$ along $f(\Sigma)$.
Let $X\subset W$ be a small neighborhood of $f(V)$.
We supply $X$ with  an ``almost normal" (to $f(V)$) foliation $\N$,
which coincides with the normal foliations to $\widehat V$ near
$f(\Sigma)$ and with the normal foliation to $f(V\setminus \Op\Sigma)$
near $f(V\setminus \Op\Sigma)$.

\subsection{Fibered wrinkled embeddings}\label{sec:fiber-wrinkl-emb}

The notion of a wrinkled embedding can be extended to the parametric case.
Considering $k$-parametric {\it families} $f_p$, $p\in B$,
of the  wrinkled embeddings we allow, in addition to
the wrinkles, their {\it embryos}.   Near each embryo $v_i\in V$
the map $f_p$ is equivalent to the map
$$Z_0(n,m):\Op_{\bbR^n} 0\to \bbR^{m}$$
given by the formula
$$(y_1,...,y_{n-1},z)\longmapsto
\left(y_1,...,y_{n-1},z^3+3|y|^2z,
\int (z^2+|y|^2)^2\,dz,0,...,0\right).\eqno (4)$$

\n
Thus, wrinkles may {\it appear and disappear}
when we consider a {\it homotopy} of wrinkled embeddings
$f_t:V\to W$.

\mn
In a more formal mode one can use, quite similar to the case of wrinkled
mappings, the ``fibered" terminology (see \cite {[EM97]}).
A {\it fibered} (over $B$) {\it map} is a commutative diagram
$$
\xymatrix{V^{k+n} \ar[rd]_p\ar[rr]^f& &W^{k+m}\ar[ld]^q\\
&B^k&}
$$

\n
where $f$ is a smooth map and $p\,$,$\,q$ are submersions. For the fibered 
map $f$ we denote by $T_BV$ and $T_BW$ the subbundles
$\Ker\,p \subset TV$ and $\Ker\,q\subset TW$.
They are tangent to foliations of $V$ and $W$ formed by
preimages $p^{-1}(b) \subset V$, $q^{-1}(b) \subset W$, $b \in B$.
We will often denote a fibered map simply by  $f:V \rightarrow W$ \,
when  $B$, $p$ and $q$ are implied from the context. Fibered homotopies, 
fibered differentials, fibered submersions, and so on
can be  naturally defined in the category of fibered maps.
For example, the {\it fibered differential} $d_Bf$ of
a fibered map $f:V \rightarrow W$ is the restriction
$$d_{B}f = df|_{T_BV} : T_{B}V \rightarrow T_{B}W\,.$$
Notice that $d_Bf$ itself has the
structure of a map fibered over $B$.

\n
Two fibered maps, $f:V\to W$ over $B$  and $f:V'\to W'$ over $B'$,
are called {\it equivalent} if there exist open subsets
$A \subset B,\,\,
A'\subset B',\,\,
Y\subset W,\,\,
Y'\subset W'$
with
$f(V)\subset Y,\,\,
p(V)\subset A,\,\,
f'(V')\subset Y',\,\,
p'(V') \subset A'$
and diffeomorphisms
$\varphi: U \rightarrow U', \; \psi : Y \rightarrow Y', \; s: A \rightarrow A'$
such that they form the following commutative diagram
$$
\xymatrix@C=14pt{V\ar[rrrrdddd]_p\ar[rrrd]^\varphi\ar[rrrrrrrr]^f&&&&&&&&Y
\ar[lllldddd]^q\ar[llld]_\psi\\&&&V' \ar[rd]_{p'}\ar[rr]^{f'}& &Y'
\ar[ld]^{q'}&&&\\&&&&A'&&&&\\&&&&&&&&\\&&&&A\ar[uu]_s&&&&}
$$
The canonical form $Z(k+n,k+m)$, being considered as a fibered map over
the space of first $k$ coordinates,
gives us the canonical form for the {\it fibered wrinkle} of the
{\it fibered wrinkled embedding}.
Thus, in a $k$-parametric family of the wrinkled embeddings
each fibered wrinkle bounds a $(k+n)$-dimensional {\it fibered}
disk, fiberwise equivalent to the standard $(k+n)$-disk in the space
$\bbR^{k+n}$, fibered over $\bbR^k$. In addition, the
canonical form $Z(k+n,k+m)$ over the half space $\bbR^{k}_-\,(y_1\leqslant 0)$
gives us the model for the fibered ``half-wrinkles'' near the boundary of $B$.

\subsection{Main theorem}\label{sec:global-integr}

The following Theorem \ref{thm:int-tang-global}, and its fibered analog
\ref{thm:int-tang-global-fib} in Section
\ref{sec:fibered} below
are the main results of the paper.
\bp
\label{thm:int-tang-global}
{\bf (Integrable approximation of a tangential rotation)}
Let $G_t:V\to \Gr_nW$ be a tangential rotation  of an embedding
$i_V:V\hookrightarrow W$.
Then there exists a homotopy of wrinkled embeddings
$f_t:V\to W$, $f_0=i_V$, such that the homotopy $\rG df_t:V\to\Gr_nW$
is arbitrarily $C^0$-close to $G_t$.
If the rotation $G_t$ is fixed on a closed subset $C\subset V$, 
then the homotopy $f_t$ can be chosen also fixed on $C$.
\ep

\n
\btl\,
{\bf Remark.} The parametric version of Theorem \ref{thm:int-tang-global}
is also true, see Section \ref{sec:fibered} below.
\btr

\mn
A small rotation $G_t:V\to \Gr_nW$
is called {\it simple}, if  $G_t(v)$ is a rotation in a $(n+1)$-dimensional
subspace $L_v\subset T_vW$ for every $v\in V$ and the angle of the rotation
$G_t(v)$ is less than $\pi/4$.
In particular, for $q=n+1$ any small rotation with the maximum
angle $<\pi/4$ is simple.
Any tangential rotation can be approximated by a finite sequence of
simple rotations and hence Theorem \ref{thm:int-tang-global}
follows from

\bp
\label{thm:int-tang-global-small}
{\bf (Integrable approximation of a simple tangential rotation)}
Let $G_t:V\to \Gr_nW$ be a simple tangential rotation of a
{\bf wrinkled} embedding $f_0:V\to W$.
Then there exists a homotopy of wrinkled embeddings
$f_t:V\to W$, such that the homotopy $\rG df_t:V\to\Gr_nW$
is arbitrarily $C^0$-close to $G_t$.
{If the rotation $G_t$ is fixed
on a closed subset $C\subset V$, then the homotopy $f_t$ can be
also chosen  fixed on $C$.}
\ep

\n
\btl\,
{\bf Remark.}
The proof of Theorem \ref{thm:int-tang-global-small}
will give us the following additional information: the homotopy $f_t$
increases the depth of $f_0$
at most by 1. In particular,
\bi
\item if $f_0$ is a smooth
embedding, then $f_t$ consists of wrinkled embeddings of depth $\leqslant 1$;
\item the depth of the final map $f_1$ in Theorem
\ref{thm:int-tang-global} is equal to the number of simple rotations
in the decomposition of the rotation $G_t$. \btr
\ei

\subsection{Local integrable approximations near wrinkles}
\label{local-integr-wrinkl}

We will distinguish the notions of {\it homotopy of wrinkled embeddings}
and their {\it isotopy}.
A homotopy $f_t$, $t\in[0,1]$, of   wrinkled embeddings
is called an {\it isotopy}, if   for all $t\in[0,1$ the wrinkled embedding $f_t$
has no embryos, i.e.  its wrinkles do not die,
and no new  wrinkles are born.
Equivalently, $f_t$ is an isotopy if $f_t=h_t\circ f_0$ where $h_t:W\to W$
is a diffeotopy.

A $C^0$-small {\it isotopy}  $f_t:V\to W$, $f_0=i_V$, of wrinkled embeddings
is called {\it graphical}, if all the images $f_t(V)$
are transversal to the normal foliation $\N$ of $V$.
As in the smooth case, the tangential rotation $G_t:V\to \Gr_nW$
of the wrinkled embedding $f_0:V\to W$ is called {\it small},
if $G_t(v)$ is transversal to $\N$ for all $t$ and $v$.
We will reformulate now   \ref{thm:int-tang-local}   for the situation,
when $f_0:V\to W$ is a {\it wrinkled} embedding and $K=\Sigma=\Sigma (f_0)$.
\bp
\label{thm:int-tang-local-wrinkled}
{\bf (Local integrable approximation of small tangential
rotation near wrinkles)}
Let $G_t:V\to \Gr_nW$ be a small tangential rotation of a
wrinkled embedding
$f_0:V\to W$.
Then there exists an arbitrarily $C^0$-small graphical isotopy
of wrinkled embeddings $f_t:V\to W$
such that the homotopy
$$\rG df_t|_{\Op \Sigma}:\Op \Sigma\to \Gr_nW$$
is arbitrarily $C^0$-close to the tangential rotation $G_t|_{\Op \Sigma}$.
\ep

\n
{\bf Proof.} The image $f_0(\Op \Sigma)\subset W$ is not a submanifold
and hence the proof of   \ref{thm:int-tang-local} formally
does not work. Let $\widehat V\subset W$ be an $n$-dimensional
submanifold, such that $f_0(\Sigma)\subset \widehat V$ and
$\widehat V$ is tangent to $f_0(V)$ along $f_0(\Sigma)$.
Let $\widehat G_t$ be an extension of the rotation $G_t|_\Sigma$
to a tangential rotation $\widehat V\to \Gr_nW$.
Apply   \ref{thm:int-tang-local} to
the pair $(\widehat V, f_0(\Sigma))$, the normal foliation $\N$
and the tangential rotation $\widehat G_t$;
denote the produced isotopy $\widehat V\to W$ by $\widehat f_t$.
Let $g_t:W\to W$ be an ambient diffeotopy for $\widehat f_t$,
such that $g_t^*\N=\N$. Then $f_t={ g}_t\circ f_0$ is the desired
isotopy of wrinkled embeddings on $\Op \Sigma$., which can be then extended
to the whole $V$   as a graphical isotopy.
\qed

\subsection{Main lemma}\label{sec:main-lemma}

Lemma \ref{thm:hypersurface}, which we prove in this section, is the main
ingredient in the proof of Theorem \ref{thm:int-tang-global-small}.

\medskip

Let us denote by $\wt\rG:S\to S^{n}$ the oriented Gaussian map  of an oriented
hypersurface  $S\subset\bbR^{n+1}$.
For the angle metric on $S^n$ denote by $U_{\e}(s)$ the open metric
$\e$-neighborhood of the north pole $s=(0,\dots,1)\in S^n$.

\mn
An oriented hypersurface $S\subset \bbR^{n+1}$ is called:
\bi
\item {\it $\e$-horizontal},
if $\,\wt\rG(S)\subset U_{\e}(s)\subset S^n$;
\item {\it graphical}, if
$\,\wt\rG(S)\subset U_{\pi/2}(s)\subset S^n$;
\item {\it $\e$-graphical}, if
$\,\wt\rG(S)\subset U_{(\pi/2)+\e}(s)\subset S^n$;
\item {\it quasi-graphical}, if
$\,\wt\rG(S)\subset U_{\pi}(s)\subset S^n$.
\ei

\mn
Similarly a wrinkled embedding $f:S\to \bbR^{n+1}$,
is called $\e$-{\it horizontal},
if its Gaussian image $\wt\rG(S)\subset S^n$ is contained  in $U_{\e}(s)$.

\mn
In what follows we will often
say {\it almost horizontal} and {\it almost graphical}
instead of  {\it $\e$-horizontal} and {\it $\e$-graphical},
assuming that $\e$ is appropriately small.

\bp
\label{thm:hypersurface}
{\bf (Approximation of an embedded hypersurface by almost horizontal
wrinkled embeddings)}
Let $S\subset \bbR^{n+1}$ be an oriented quasi-graphical
hypersurface, such that $S$
is almost horizontal near the boundary $\p S$.
Then there exists a $C^0$-approximation of the embedding
$i_S:S\hookrightarrow \bbR^{n+1}$
by an almost horizontal wrinkled embedding
$f:S\to \bbR^{n+1}$ of depth 1, such that
$f$ coincides with $i_S$ near $\partial S$.
\ep

\n
{\bf Proof.} Cutting from $S$ a small neighborhood of the critical
set of the function $h=x_{n+1}:S\to\bbR$, where $S$ is already almost
horizontal, we may  consider this new $S$,
and thus assume that the function $h$ has no critical points on $S$.
Let $\wt S\subset \Int S$ be a slightly smaller compact  hypersurface
such that $S\setminus\wt S$ is still almost horizontal.
On a neighborhood $\Op S$ consider a
non-vanishing vector field $v$ which
is transversal to  $S$, defines its given co-orientation, and horizontal,
i.e. tangent to the level sets of the function $h$.
Let $v^t$ denotes the flow of $v$. We may assume that for $z\in S$ the
flow-line  $v^t(z)$ is  defined for all $|t|\leq 1$.

\mn
Consider two foliations on $S$: the $1$-dimensional
foliation $\G$, formed by the gradient trajectories of the function
$h$, and  the $(n-1)$-dimensional foliation $\L$, formed by
level surfaces of the function $h$.

\mn
A subset $C\subset S$ is called {\it cylindrical},
if $C$ equipped with $\L$ and $\G$ is diffeomorphic to $D^{n-1}\times D^1$
equipped with the standard horizontal and vertical foliations. We call $C$
{\it special} if, in addition, $C$ is almost horizontal near its top and bottom.

\mn
Fix a covering of $\wt S$ by special cylindrical sets $C_\alpha\subset S$,
$\alpha =1,...,K$, and fix the respective parameterizations
$\varphi_\alpha:D^{n-1}\times D^1\to  C_\alpha$,
which send the positive direction on $D^1$   to the
directions on $\G$   defined by the gradient  vector field $\nabla h$.
Given a  function $\tau:S\to\bbR$, $|\tau|\leq 1$, we will denote by
$I^\tau_S$ a perturbation of the inclusion $i_S$ defined by the formula
$$I^\tau_S(z)=v^{\tau(z)}(z),\;\; z\in S.$$

\begin{figure}[hi]
\centerline{\psfig{figure=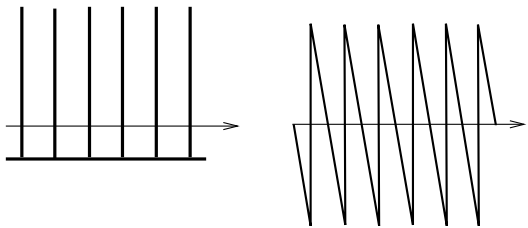,height=30mm}}
\caption {\small Functions $w$ and $W$}
\label{we08}
\end{figure}

\begin{figure}[hi]
\centerline{\psfig{figure=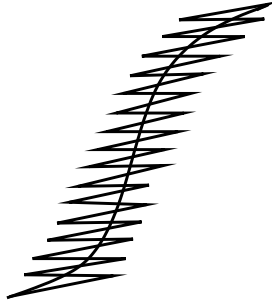,height=40mm}}
\caption {\small Local wrinkling, the case $n=1$}
\label{we09}
\end{figure}

\n
To clarify the geometric meaning of the construction of the necessary
approximation
we will give first a  slightly imprecise description without technical details.
Let $x_i\in D^1$, $i=1,...,N$, be a finite set of points such that $x_1=-1+d/2$,
$x_N=1-d/2$ and $x_{i+1}-x_i=d$. Set $w(x)=\sum_i\d(x-x_i)-1/d$, where
$\d(x)$ is  the  $\d$-function, and
define $W(x):=\int_0^xw(x)dx$ (see Fig.\ref{we08}).
Take a cut-off function $\lambda:D^{n-1}\times D^1 \to\bbR_+$,
equal to $1$ on a slightly lesser than $D^{n-1}\times D^1$
subset $A\subset D^{n-1}\times D^1$
and equal to $0$ on $\Op \p (D^{n-1}\times D^1)$, and define
a function $\tau:D^{n-1}\times D^1 \to\bbR$ by the formula
$\tau(r,x)=\lambda(r,x)W(x)$ and a function
$\tau_\alpha:C_\alpha\to\bbR^1$ by the formula
$\tau_{\alpha}(c)=\gamma(\tau\circ\varphi_\alpha^{-1}(c))$, $\gamma\in \bbR_+$.
For sufficiently small $d$ and $\gamma$
the map $I^{\tau_\alpha}_S$ (see Fig.\ref{we09}) is a good candidate
for the wrinkled embedding.
Take this candidate for a moment; then
assuming that all the levels $h=\varphi_\alpha(x_i)$ for all $C_\alpha$
are distinct from each other, we can apply the described
local goffering  {\it simultaneously} for all $C_\alpha$ and thus get the desired
global approximation of $S$.
Unfortunately, our map $I^{\tau_\alpha}_S$ does not match the canonical form (1)
and requires an additional perturbation near the singularity. Rather than doing
this, we describe below the whole construction  again  in a slightly modified
form with all the details.
%To construct the required approximation by {\it genuine} wrinkled embeddings let
%us introduce some special models.

\begin{figure}[hi]
\centerline{\psfig{figure=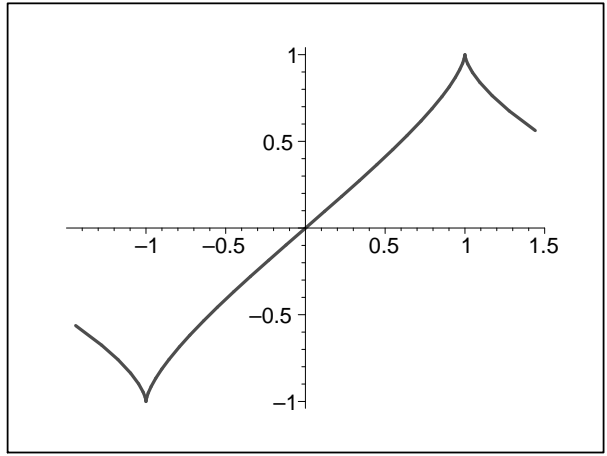,height=60mm}}
\caption {\small Function $a_t$, $t=1$}
\label{m1}
\end{figure}

\begin{figure}[hi]
\centerline{\psfig{figure=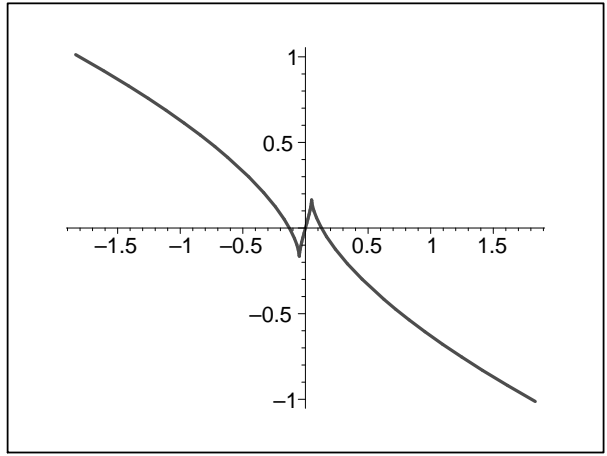,height=60mm}}
\caption {\small Function $a_t$, $t=0.3$}
\label{m2}
\end{figure}

\mn
First, choose a family of  curves
$A_t\subset \bbR^2$, $t\in\bbR$, which is given by parametric equations
$$
x(u,t)=\frac{15}8\int\limits_0^u(u^2-t)^2 du\,,\;\;\;
y(u,t)=-\frac12(u^3-3tu).
$$
The curve $A_t$ is a graph of a continuous function $a_t:\bbR\to\bbR$ which
is smooth for $t<0$ and smooth on $\bbR\setminus\{-\sqrt t,\sqrt t\}$
for $t\geq 0$, see Fig.\,\ref{m1},\ref{m2}. Note that 
\bi
\item [(a)] $a_1(\pm 1)=\pm 1$;
\item [(b)] the composition $a_t(x(u,t))$ is a smooth function.
\ei
\begin{figure}[hi]
\centerline{\psfig{figure=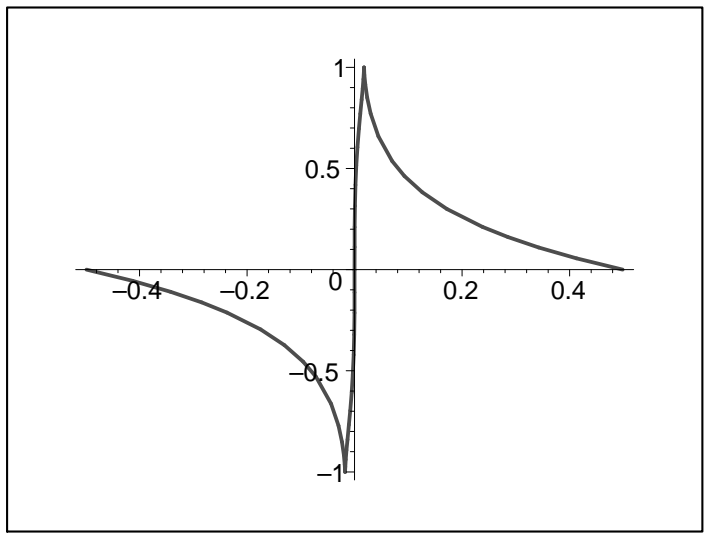,height=70mm}}
\caption {\small Function $\W_{\sigma,1}$,\,\,\,$x\in {[-\frac12,\frac12]}$,\,
$\sigma\simeq 0.02$}
\label{m3}
\end{figure}
Choose $\sigma\in(0,\frac18)$  and consider a family of odd $1$-periodic
functions $\W_{\sigma,\,t}:\bbR\to\bbR$ with the following properties:
$$\W_{\sigma,\,t}(x)=
 \begin{cases} a_t\left({\frac x\sigma}\right)&\hbox{ for}\;\;
   x\in\Op[-\sigma,\sigma]\cr 0&\hbox{ for}\;\;  x=\frac12\cr
 \end{cases}
$$
$$
 \frac{d\W_{\sigma,1}}{dx}(x)
 \begin{cases}
 \geqslant 3&\hbox{ for} \;\;x\in(-\sigma,\sigma)\cr
 \in[-2,0)&\hbox{ for} \;\;x\in[2\sigma, \frac12-2\sigma]\cr
\leqslant -2&\hbox{ for}\;\; x\in[\sigma, 2\sigma]\cr
\end{cases}
$$
(see Fig.\,\ref{m3}). Note that the composition 
$\W_{\sigma,\,t}(a_t(x(u,t))$ is a smooth function.
Pick an $\e>0$ and integer $N>0$ such that the cylinders
$$\wt C_\alpha=\varphi_\alpha(D^{n-1}_{1-2\e}\times D^1_{1-1/N})$$
are still special and cover the surface $\wt S$. We will keep $\e$ fixed while
$N$ will need to be increased to achieve the required approximation.
Choose a non-strictly decreasing smooth  function $\beta:[0,1]\to \bbR$ 
such that
$$
 \beta(r)=
 \begin{cases}  1&   r\in [0,1-3\e];\cr
1- r-\e& r\in[1-2\e,1].\cr
 \end{cases}
$$
Given an integer $N>0$, pick also a cut-off function $\lambda:[0,1]\to[0,1]$
which is equal to $1$ on $[1-\rho]$ and to $0$ near $1$, 
where $\rho=\min(\frac{1}{2N},\,\frac{\e}{2})$,  and
define a function $$\tau=\tau_{N,\sigma}:D^{n-1}\times D^1\to\bbR$$
by the formula
$$\tau(r,x)=\gamma\lambda(|x|)\lambda(|r|)
\W_{\sigma,\,\beta(|r|)}\left(\frac{2N+1}2x\right),\;\;
(r,x)\in D^{n-1}\times D^1,\,\,\,\gamma\in \bbR_+. $$
Let us push-forward the function $\tau $ to $S$ by $\varphi_\alpha$,
$$\tau_{\alpha}:=\tau\circ\varphi_\alpha^{-1}.$$
\bp
\label{lm:model-wrinkling}
{\bf (Local wrinkling)}
For sufficiently big $N$ and sufficiently small $\gamma$
the map $I^{\tau_\alpha}_S\circ g:S\to \bbR^{n+1}$ is a wrinkled embedding.
Here $g:S\to S$ is an appropriate smooth homeomorphism,
which makes the composition smooth.
\ep
{\bf Remark.} Of course, $I^{\tau_\alpha}_S(S)=(I^{\tau_\alpha}_S\circ g)(S)$.
However, the map $I^{\tau_\alpha}_S$ itself is not smooth.
Using the property (b) of the function $a_t$ 
one can write down the reparameterization $g$ explicitly.

\mn 
Let $  D_{\bbQ}^1$ be the set of rational points in $D^1$.
We can choose parameterizations
$\varphi_\alpha$ is such a way that the images 
$\varphi_\alpha(D^{n-1}\times D^1_{\bbQ})\subset S$
are pairwise disjoint. Then for a sufficiently large $N$ the images
$$\varphi_\alpha(\bigcup\limits_{-N}^N[x_k-\wt\sigma,x_k+\wt\sigma])
\subset S,$$
where
$ x_k=\frac {2k}{2N+1}, \wt\sigma=\frac {4\sigma}{2N+1}$,
are pairwise disjoint as well.

\n
Choose a partition of unity $\sum\limits_1^K\eta_\alpha=1$  on $\wt S$ realized
by functions $\eta_\alpha$ supported in
$\Int\, C_\alpha$ and such that $\eta_\alpha|_{\wt C_\alpha}>0$.
Finally, for  a sufficiently large $N$ and sufficiently small $\gamma>0$,
the map $I_S^\Psi\circ g$, where
$$\Psi :=\sum\limits_{\alpha=1}^K\eta_\alpha  {\tau_{\alpha}},$$
is the required approximation of $i_S$ by an almost horizontal wrinkled
embedding. Here, like in \ref{lm:model-wrinkling}, $g:S\to S$ is an appropriate 
reparameterization.  
\qed

\mn
We will need the following {\it parametric} version of 
Lemma \ref{thm:hypersurface}:

\bp
\label{thm:hypersurface-param}
{\bf (Approximation of  a family of embedded hypersurfaces
by a family of almost horizontal
wrinkled embeddings)}
Let $S_t\subset \bbR^{n+1}$, $t\in I$, be a family of
oriented quasi-graphical
hypersurfaces, such that $S_t$ is almost horizontal for
$t=0$ and $S_t$
is almost horizontal near the boundary $\p S_t$ for all
$t\in I$.
Then there exists a $C^0$-approximation of the family of
embeddings $i_{S_t}:S_t\hookrightarrow \bbR^{n+1}$
by a family of almost horizontal wrinkled embeddings
$f_t:S_t\to \bbR^{n+1}$ of depth $\leqslant 1$, such that
$f_t$ coincide with $i_{S_t}$ for $t=0$ and
$f_t$ coincide with $i_{S_t}$ near $\partial S$ for all $t\in I$.
\ep

\n
In order to prove this version, we need to do the following
modification in the previous proof. The family
$S_t$, $t\in I$, can be considered as a fibered over $I$
quasi-graphical hypersurface $S$ in $I\times \bbR^{n+1}$. Then,
all above  constructions should be done in the  fibered category. In particular,
instead of the cylinder $D^n\times D^1$ we need to use
the fibered cylinder $I\times D^{n-1}\times D^1$
for the parameterization of the fibered cylindrical sets.

The hypersurface $S$ is almost
horizontal everywhere near the boundary $\p S$ except
the right side, $S_1\subset \p S$. Hence, in order to use the above scheme,
we will work with  the doubled, fibered over $[0,2]$
family $S_t, t\in[0,2]], $ where $S_t=S_{2-t}$ for $t\in[1,2]$, while
making all the above  constructions equivariant with respect to the involution
$t\mapsto 2-t$, $t\in[0,2]$.

\subsection{Proof of Theorem \ref{thm:int-tang-global-small}}
\label{sec:proof-int}

In order to get rid of non-essential details, we will
consider only the case $(W,g_W)=(\bbR^q,dx_1^2+...+dx_m^2)$.
The proof can be easily adjusted to the general case.

\mn
Working with a triangulation $\Delta$ of the manifold $V$,
we always assume that $S'\subset S=\Sigma(f_0)$ is contained in the
$(n-2)$-skeleton of $\Delta$ and $S$ in contained in its
$(n-1)$-skeleton.
Given a triangulation $\Delta$ of $V$ and a map $G:V\to \Gr_n\bbR^m$,
we will denote by $G^\Delta$ the
piecewise constant map $G^\Delta:V\to \Gr_n\bbR^m$, defined on each
$n$-simplex $\Delta_i ^n$ by the condition that
$G^\Delta(v)$ is parallel to $G(v_i)$, where $v_i$ is the
barycenter of the simplex $\Delta^n_i$.
The map $G^\Delta$ is multivalued over the $(n-1)$-skeleton of $\Delta$.

\mn
{\bf The case $\bf m=n+1$.}
We can choose a triangulation $\Delta$ of $V$
with sufficiently small simplices ({\it fine} triangulation),
such that for any $i$ the image $f_0(\Delta^n_i)$ is arbitrarily
$C^1$-close to $G_0(v_i)$ and for any $t$ the map $G_t^\Delta$
is arbitrarily $C^0$-close to $G_t$.
In what follows we will approximate the homotopy $G_t^\Delta$ instead
of $G_t$, and therefore the approximation $G_t\approx G_t^\Delta$ must
have at least the same order as the desired approximation
$G_t\approx \rG df_t$.

\mn
First, we apply Theorems \ref{thm:int-tang-local} and
\ref{thm:int-tang-local-wrinkled}
to the $(n-1)$-skeleton $\Delta^{n-1}$ and construct a graphical
isotopy $\wt f_t:V\to \bbR^{n+1}$, $\wt f_0=f_0$,
such that $\rG d\wt f_t|_{\Op \Delta^{n-1}}$
is arbitrarily $C^0$-close to $G_t|_{\Op\Delta^{n-1}}$.
For every $i$ and $t$
the hypersurface $\wt f_t(\Delta^n_i)$,
i.e. the image of $f_0(\Delta^n_i)$ under the graphical isotopy,
is almost graphical
with respect to $\rG_0(v_i)$, and hence it is  quasi-graphical
with respect to the  hyperplane $G_t(v_i)$ because the angle
between $\rG_0(v_i)$ and $G_t(v_i)$ is less then $\pi/4$.
Then, $\wt f_t(\Delta^n_i)$ is almost horizontal with respect
to $G_t(v_i)$ near the boundary $\p\wt f_t(\Delta^n_i)$.
Therefore, over each $n$-simplex $\Delta^n_i$ we can apply
Lemma \ref{thm:hypersurface-param} with $S_t=\wt f_t(\Delta^n_i)$
and $G_t(v_i)$ as the horizontal hyperplane,
i.e. for every $t$ the hyperplane
$G_t(v_i)$ plays the role of the horizontal hyperplane
$\bbR^n\subset \bbR^{n+1}$.
Thus we can deform the isotopy of embeddings
$\wt f_t$ to the desired homotopy $f_t$ of the wrinkled embeddings.

\mn
{\bf The case $\bf m>n+1$.}
Let $\Delta$ be a fine triangulation of $V$.
For every $i$ the image $f_0(\Delta^n_i)$ is arbitrarily
$C^1$-close to $G_0(v_i)$ and, therefore,  we can work over each simplex
$\Delta_i^n$ with the projection of the isotopy
$\wt f_t$ to $L_{v_i}\simeq \bbR^{n+1}$ (the $(n+1)$-dimensional
subspace where the rotation $G_t$ goes on) exactly as in the
previous case, keeping the coordinates in $L_{v_i}^\bot$
unchanged.
\qed

\subsection{Fibered case}\label{sec:fibered}

Let us formulate and discuss the parametric version of Theorem
\ref{thm:int-tang-global}.
For a fibration $q:W\to B$ denote by $\Gr_nW_B$ the Grassmannian
of $n$-planes tangent to the fibers of the fibration $q$.

\bp
\label{thm:int-tang-global-fib}
{\bf (Global integrable approximation of fibered tangential rotation)}
Let $G_t:V\to \Gr_nW$ be a fibered tangential rotation  of a fibered (over $B$)
embedding $i_V:V\hookrightarrow W$.
Then there exists a homotopy of fibered wrinkled embeddings
$f_t:V\to W$, $f_0=i_V$, such that the (fibered) homotopy
$\rG d_Bf_t:V\to\Gr_nW_B$
is arbitrarily $C^0$-close to $G_t$. If the rotation $G_t$ is fixed
on a closed subset $C\subset V$, then the homotopy $f_t$ can be chosen
also fixed on $C$.
In particular, if $G_t$ is fixed over a closed subset $B'\subset B$,
then the homotopy $f_t$ can be chosen also fixed over $B'$.
\ep

\n
This theorem follows from

\bp
\label{thm:int-tang-global-small-fib}
{\bf (Global integrable approximation of simple fibered tangential rotation)}
Let $G_t:V\to \Gr_nW$ be a simple fibered tangential rotation of a
fibered {\bf wrinkled} embedding $f_0:V\to W$.
Then there exists a homotopy of fibered wrinkled embeddings
$f_t:V\to W$, such that the (fibered) homotopy $\rG d_Bf_t:V\to\Gr_nW_B$
is arbitrarily $C^0$-close to $G_t$.
\ep

\n
Now we need a version of the main lemma \ref{thm:hypersurface}
which is parametric with respect to  both,  time and space:

\bp
\label{thm:hypersurface-param-fib}
{\bf (Approximation of family of fibered embedded hypersurfaces
by family of almost horizontal fibered wrinkled embeddings)}
Let $S_t\subset B^k\times \bbR^{n+1}$, $t\in I$, be a family of
fibered over $B^k$ oriented quasi-graphical
hypersurfaces, such that $S_t$ is almost horizontal for
$t=0$ and $S_t$ is almost horizontal near the boundary $\p S_t$
and over $\p B$ for all $t\in I$.
Then there exists a $C^0$-approximation of the family of fibered
embeddings $i_{S_t}:S_t\hookrightarrow B\times\bbR^{n+1}$
by a family of almost horizontal fibered wrinkled embeddings
$f_t:S_t\to \bbR^{n+1}$ of depth $\leqslant 1$, such that
$f_t$ coincides with $i_{S_t}$ for $t=0$ and
$f_t$ coincides with $i_{S_t}$ near $\partial S$ and over $\p B$
for all $t\in I$.
\ep

\n
The proof of the time-parameric version \ref{thm:hypersurface-param}
of the main lemma can be rewritten for this fibered case.
In the proof of Theorem \ref{thm:int-tang-global-small-fib}, in order to
apply the fibered version of the main lemma, the fine triangulation $\Delta$
ought to be transversal to the fibers on $\Int\, \Delta_i$ for all its
simplices $\Delta_i\subset S\setminus \Op \p S$. 
Let us point that the existence of such a triangulation is non-obvious
and was proven in a more general case (for foliation) by W.~Thurston 
in \cite{[Th74]}.

\subsection{Double folds}\label{sec:double}

A (smooth) map $f:V^n\to W^m\,,\,\,n<m\,,$
is called {\it folded embedding}, if
\bi
\item $f$ is a topological embedding;
\item the singularity $\Sigma=\Sigma_f$ of the map $f$ is an $(n-1)$-dimensional
submanifold in $V$;
\item near each connected components $S$ of $\Sigma$ the map $f$ is
equivalent to the map
$$\Op_{{S}\times\bbR^1}S\times 0\to S\times\bbR^{m-n+1}$$
given by the formula
$(y,z)\to (y,z^2, z^3, 0,...,0)$, where $y\in S$.
\ei
Thus,  folded embeddings have only {\it two-fold} corners,
ike on Fig.\ref{we01}. The submanifolds $S$ will be called the {\it folds}
of the folded embedding $f$.

We say that a folded embedding $f$ has  {\it spherical double-folds} if all its
folds diffeomorphic to the the $(n-1)$-sphere,
and  organized in pairs $(S_0,S_1)$
which bounds annuli diffeomorphic  to $S^{n-1}\times I$.
These annuli are allowed to be {\it nested}, i.e. the annulus bounded
by one double fold may contain the annulus associated with another double fold.

Considering families of folded embeddings with spherical double folds we will
also allow, similar to the case of wrinkled embeddings, {\it embryo}
double folds, so that a double fold could die or be born during the deformation.
The local model  for a  map
$\Op_{{S}\times\bbR^1}S\times 0\to S\times\bbR^{m-n+1}$
near an embryo double fold is given by the formula
$$(y ,z)\longmapsto
(y ,z^3,z^5,0,...,0). $$

Theorem \ref{thm:int-tang-global} has an equivalent reformulation
for folded embeddings with spherical  double folds.
%%%%%%%%%%%%%%%%
\bp
\label{thm:int-tang-global-folds}
{\bf (Global integrable approximation of tangential rotation by folded embeddings
with spherical double folds)}
Let $G_t:V\to \Gr_nW$ be a tangential rotation  of an embedding
$i_V:V\hookrightarrow W$.
Then there exists a homotopy of folded embeddings with spherical double folds
$f_t:V\to W$, $f_0=i_V$, such that the homotopy $\rG df_t:V\to\Gr_nW$
is arbitrarily $C^0$-close to $G_t$.
If the rotation $G_t$ is fixed
on a closed subset $C\subset V$, then the homotopy $f_t$ can be chosen
also fixed on $C$.
\ep
%%%%%%%%%%%%%%%%%%%
A fibered version of this theorem also holds.

\begin{figure}[hi]
\centerline {\psfig{figure=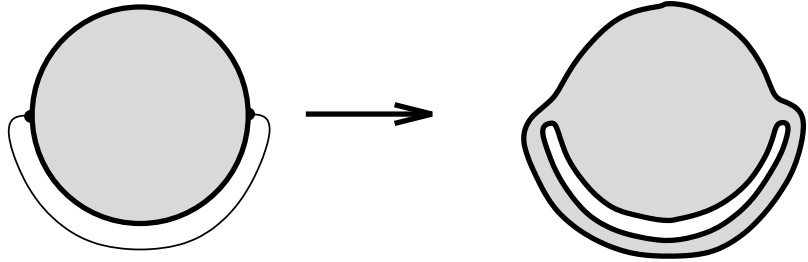,height=35mm}}
\caption {\small Whitney surgery (in the pre-image); $n=2$}
\label{we13}
\end{figure}

\mn The following surgery construction allows one to deduce
\ref{thm:int-tang-global-folds} from \ref{thm:int-tang-global}.
Let $f:V^n\to W^m$ be a wrinkled embedding, $n>1$.
One can modify each wrinkle of $f$ by a connected sum construction
for cusps (Whitney surgery) such that the resulting
map $\wt f:V^n\to W^m$ will be a folded embedding with
{\it spherical double folds} $S^{n-1}\times S^0$.
For maps $V^2\to W^m$ the construction (in the pre-image)
is shown on Fig.\ref{we13}.
Next two propositions contain  a formal description of the construction,
see \cite{[El72]} or \cite{[EM98]} for more details.

\bp
\label{thm:preparation}
{\bf (Preparation for the cusp surgery)}
Let
$$Z(n,q): \Op_{\bbR^n} S^{n-1}\to \bbR^{m}$$
be the standard wrinkled embedding with the wrinkle $S^{n-1}\subset \bbR^n$.
If $n>1$ then there exists an embedding
$h:D^{n-1}\to \Op_{\bbR^n} S^{n-1}$ such that

\begin{itemize}
\item $h(\partial D^{n-1})=S^{n-2}$, the cusp (three-fold points) of the wrinkle;
\item $h(D^{n-1}\setminus \Int\, D^{n-1}_{1-\delta})=
(D^{n-1}_{1+\delta}\setminus\Int\, D^{n-1})\times 0\subset
\bbR^{n-1}\times\bbR^1$;
\item $h(\Int\, D^{n-1})$ does not intersect the wrinkle $S^{n-1}$.
\end{itemize}
\ep
\bp
\label{thm:surgery}
{ (\bf Surgery of cusps)}
Let $h:D^{n-1}\to \Op_{\bbR^n} S^{n-1}$ be an embedding, as in
\ref{thm:preparation}.
There exists a $C^0$-small perturbation of the map
$Z(n,m)$ in an arbitrarily small
neighborhood of the embedded disk $h(D^{n-1})$ such that
the resulting map $\widetilde Z(n,m)$
is a folded embedding with two spherical folds.
\ep

\n
The construction also implies that each double fold $S^{n-1}\times S^0$
of the folded embeddings $\wt f$ bounds an annulus $S^{n-1}\times D^1$ in
$V$.

\section{Applications}\label{appl}

\subsection[Homotopy principle for directed wrinkled embeddings]
{Homotopy principle for directed wrinkled \\
embeddings}\label{sec:direct}

Using the wrinkled embeddings one can reformulate the $h$-principle
for $A$-directed embeddings of {\it open} manifolds (see 4.5.1 in \cite{[EM02]})
to the case of $A$-directed  wrinkled embeddings of {\it closed} manifolds.

\bp
\label{thm:dir-wrinkled}
{\bf (A-directed wrinkled embeddings of closed manifolds)}
If $\;A\subset \Gr_nW\;$ is an open subset and $f_0:V\hookrightarrow W$
is an embedding whose tangential lift
$$\;G_0=\rG df_0:V\to \Gr_nW$$
is homotopic to a map
$$G_1:V\to A\subset \Gr_nW\,,$$
then there exists a homotopy of wrinkled embeddings
$f_t:V\to W$ such that $f_1:V\to W$ is
an $A$-directed wrinkled embedding.
Such a homotopy can be
chosen arbitrarily $C^0$-close to $f_0$.
\ep

\n
One can formulate \ref{thm:dir-wrinkled} also in the case when
$A$ is an open subset in the Grassmannian of {\it oriented}
$n$-planes in $W$. In both cases the h-principle holds also in the
relative and parametric versions. Then, one can reformulate this
$h$-principle for the folded embeddings with spherical double folds.

\subsection{Embeddings into foliations}
\label{sec:embed-to-fol}

{\bf A. Wrinkled mappings and generalized wrinkled mappings.}
In  Section \ref{appendix} below we provide for a convenience of the
reader the basic definitions and results of the wrinkling theory from
\cite{[EM97]} and \cite{[EM98]}.  However, in this section we talk about
{\it generalized} wrinkled maps.
Let us explain here the difference between the two notions.

\begin{figure}[hi]
\centerline{\psfig{figure=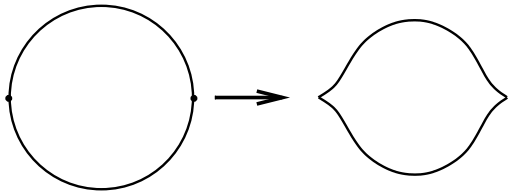,height=25mm}}
\caption {\small Wrinkle and its image; $n=q=2$}
\label{we04}
\end{figure}

\mn
A wrinkled mapping $f: V^n\to W^q\,,\,\,\, n\geqslant q\,$,
by definition, is a map whose singularity set consists of {\it $(n,q)$-wrinkles}
(= {\it wrinkles}), where each wrinkle is a $(q-1)$-dimensional
sphere $S_i\subset V$ of {\it fold} points divided by an equator
of {\it cusp} points, see Fig.\,\ref{we04}. 
In addition, each wrinkle $S_i\subset V$ is required to bound an embedded
$q$-dimensional disk $D_i\subset V$,
such that {\it the restriction of the map $f$
to $D_i\setminus S_i$ is an equidimensional embedding}.
\begin{figure}[hi]
\centerline{\psfig{figure=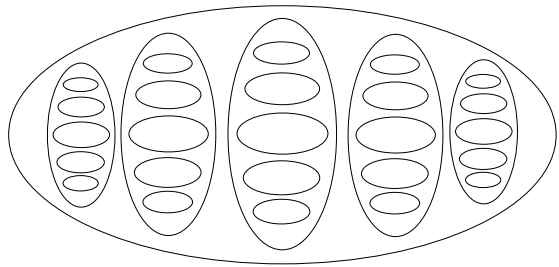,height=35mm}}
\caption {\small Nested wrinkles in the pre-image (i.e. in $V$); $n=2$}
\label{we05}
\end{figure}

\mn
In the case of  a {\it  generalized wrinkled mapping} 
$f$ each $(n,q)$-wrinkle $S_i$
also bounds a disk $D_i$ in $V$.
However, the restriction $f|_{D_i}$  {\it is not required to be a
smooth embedding on $D_i\setminus S_i$}. Instead, this restriction is allowed
itself to be wrinkled, see Fig.\ref{we05}. See Section 
\ref{sec:generalized wrinkles} below for more details. 

\mn
{\bf B. Main theorem (statement).}
Let $(W^m,\F)$ be a foliated manifold, $\codim\,\F=q\leqslant n$.
We say that the singularities of an embedding 
$f:V^n\to (W,\F)$ with respect to foliation $\F$ 
are {\it generalized wrinkles},
if $f$ is a {\it generalized wrinkled map with respect to $\F$},
see Section \ref{sec:generalized wrinkles}.

\bp
\label{thm:embedding-foliation}
{\bf (Embeddings into foliations)}
Let $\F_p$, $p\in B$, be a family of foliations on $W^m$,
${\rm codim}\,\F_p=q\leqslant n$,
and $f_p:V^n\to (W^m,\F_p)$ be a family of embeddings
such that $f_p$ is transversal to $\F_p$ for $p\in \Op \p B$.
Suppose, in addition,  that there exists a family of 
tangential rotations $G^t_p:V\to \Gr_nW$, $p\in B$, $t\in I$, such that
\bi
\item $G^t_p$ is constant for $p\in \Op \p B$;
\item $G^1_p$ is transversal to the foliation $\F_p$ for all $p$.
\ei
Then there exists a family of embeddings $f^t_p$, $p\in B$, $t\in I$,
such that
\bi
\item $f^0_p=f_p$ for all $p$ and $f^t_p=f_p$ for $p\in \Op \p B$;
\item for any $p$ the singularities of the $f^1_p$ with respect to
$\F_p$ are generalized $(n,q)$-wrinkles and embryos.
\ei
\ep

\btl\,
{\bf Remarks.}
{\bf 1.} The family $f^t_p$ is $C^0$-close to the
constant family $\bar f^t_p=f^0_p$.

\n
{\bf 2.} The theorem is true also relative to a closed subset
$A\subset V$, i.e. in the situation when the embeddings $f_p$
are already transversal to $\F$ on $\Op_V A$.

\n
{\bf 3.} In a more formal way, one can say that $f_p^1$ is 
a {\it fibered generalized wrinkled map with respect to $\F_p$\,}.
\btr

\mn
\btl\,
{\bf Examples.}
{\bf 1.} {\it Let $f_p:S^{n-1}\times I\to \bbR^{n+1}$, $p\in D^k$,
be a family of embeddings, such that $f_p$ is the standard inclusion
$i_{S^{n-1}\times I}:S^{n-1}\times I\hookrightarrow\bbR^{n+1}$
for $p\in \Op \p D^k$
and $f_p=i_{S^{n-1}\times I}$ near $\p (S^{n-1}\times I)$ for all $p$.
Then there exists an isotopy of the family $f_p$ which is fixed near
$\p S^{n-1}\times I$ and for $p\in\p D^k$ and  such that
the projections of the resulting cylinders to the axis $x_{n+1}$
have only Morse and birth-death
type singularities.}

Indeed, the  corresponding homotopical condition  for the family of normal
vector fields is automotically holds here, as it was observed by A. Douady
and F. Laudenbach, see \cite{[La76]} and \cite{[EM00]}. Similarly, we have

\n
{\bf 2.} {\it Let $f_p:D^n\to \bbR^{n+1}$, $p\in D^k$ be a family of
embeddings, such that $f_p$ is the standard inclusion
$i_{D^n}:D^n\hookrightarrow\bbR^{n+1}$
for $p\in \Op \p D^k$, and $f_p=i_{D^n}$ near $\p D^n$ for all $p$.
%Suppose that the family of normal vector fields to the images of $D^n$
%is homotopic to a family of non-horizontal vector fields.
Then there exists an isotopy of the family $f_p$, which is fixed near $\p D^n$
and for $p\in\p D^k,$   such that
all the singularities of the projections of the resulting disks on $\bbR^n$
are generalized wrinkles and embryos.}
\btr

\mn
 Theorem \ref{thm:embedding-foliation} is proved below   in 
{\bf D}.

\mn
{\bf C. $\F$-regularization (lemmas).} 
For the canonical form 
$$Z(n,n+1):\Op_{\bbR^n} S^{n-1}\to \bbR^{n+1}$$
of the wrinkled embedding (see \ref{sec:wr-embed-formal}\,A)
the {\it regularizing foliation} is, by definition, the 
one-dimensional affine foliation in $\bbR^{n+1}$,
parallel to the axis $x_{n+1}$.

\mn
In order to make the situation more transparent and simplify the notation,
we formulate our lemmas  only in the case  $m=n+1$.
All statements remain true also for an arbitrary $m\geq n+1$.

\bp
\label{thm:F-regularization-1}
Let $f:V^n\to (W^{n+1},\F)$ be a wrinkled embedding, transversal to
the foliation $\F$, \,$\dim\,\F=1$ ($\codim\,\F=n$).
Then there exists a regularization $\wt f$ of the wrinkled embedding 
$f$ such that the singularities of the smooth embedding
$\wt f$ with respect to $\F$ are {\it generalized wrinkles}.
\ep
{\bf Proof.} We can chose the canonical coordinates near each wrinkle
$f(S_i)$ such that $\F$ will be the regularizing foliation
and then apply the standard regularization
(see \ref{sec:wr-embed-formal}\,B).
\qed

\mn
Let $f:V^n\to (W^{n+1},\F)$ be a wrinkled embedding.
Denote by $\widehat V_i$ an $n$-dimensional submanifold in $W^m$,
such that $f(S_i)\subset \widehat V_i$ and
$\widehat V_i$ is tangent to $f(\Op S_i)$ along $f(S_i)$
(see \ref{sec:wr-embed-formal}\,D).
Notice that if $f$ is transversal to $\F$ then $\widehat V_i$
is transversal to $\F$ near $f(S_i)$.  

\bp
\label{thm:F-regularization-2}
Let $f:V^n\to (W^{n+1},\F)$ be a wrinkled embedding, transversal to
the foliation $\F$, $\dim\,\F\geqslant 2$ ($\codim\,\F\leqslant n-1$). 
Suppose, in addition,  that  the restriction 
$f|_{S_i}:S_i\to \widehat V_i$ is transversal to 
the foliation $\F\cap \widehat V_i$ for all $i$.
Then any regularization $\wt f$ of the wrinkled embedding $f$
(i.e. for any choice of the canonical coordinates) 
gives us a smooth embedding transversal to $\F$.\qed
\ep

\begin{figure}[hi]
\centerline {\psfig{figure=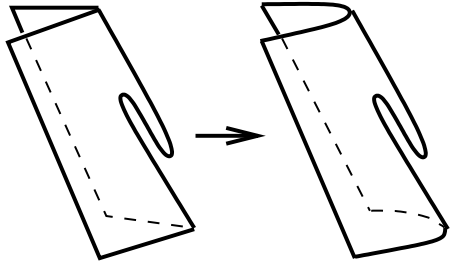,height=35mm}}
\caption {\small Regularization near the $\,(1,1)$-wrinkle}
\label{we12}
\end{figure}

\bp
\label{thm:F-regularization-3}
Let $f:V^n\to (W^{n+1},\F)$ be a wrinkled embedding, transversal to
the foliation $\F$, $\dim\,\F\geqslant 2$ ($\codim\,\F\leqslant n-1$).
Suppose, in addition, that  the singularities of the restriction
$f|_{S_i}:S_i\to \widehat V_i$ with respect to $\F\cap \widehat V_i$
are generalized $\,(n-1,q)$-wrinkles for all $i$. Then there exists a regularization 
$\wt f$ of the map $f$ such that the singularities of the smooth embedding
$\wt f$ with respect to $\F$ are generalized $\,(n,q)$-wrinkles.
\ep
{\bf Proof.} We can choose the canonical coordinates near each wrinkle
$f(S_i)$ in such a way that the one-dimensional regularizing foliation  
is inscribed into $\F$, and then apply the standard regularization.
Such a regularization adds $\,\pm t^2$ to the canonical form of any 
generalized $\,(n-1,q)$-wrinkle and thus transform 
it to a $\,(n,q)$-wrinkle. See Fig.\,\ref{we12}.
\qed

\n
{\bf Remark.}  Lemmas \ref{thm:F-regularization-1}-
\ref{thm:F-regularization-3} remain  true in the parametric form.

\mn
{\bf D. Main theorem (proof).}
We consider only a non-parametric
situation, i.e. when $B$ is just a point. The proof in the parametric
case is similar with a systematic use of the fibered terminology.
Thus, we will drop ``$p$" from  the notation. Moreover, 
we will consider only the case $m=n+1$; the proof can be easily 
rewritten for any $m>n$. 

\mn
First of all, let us apply Theorem \ref{thm:int-tang-global} to the family
$G^t:V^n\to\Gr_nW^{n+1}$ and construct a family of wrinkled
embeddings $\widehat f^t:V^n\to W^{n+1}$, $t\in [0,1]$, transversal to $\F$.

\mn
a) {\sl The base of  the induction:} $\dim\,\F=1\,$  (the equidimensional case).
We regularize the family $\widehat f^t$ in such a way that for $\widehat f^1$
our regularization is the $\F$-regularization as in Lemma 
\ref{thm:F-regularization-1}, and get the required family $f^t$
of embeddings where $f^1$ is a generalized wrinkled map with respect to $\F$.

\mn
b)  {\sl Wrinkling of wrinkles}: $\dim\,\F>1\,$.
Let $\widehat V_i\subset W$ be an $n$-dimensional submanifold, such
that $\Sigma_i=\widehat f^1(S_i)\subset \widehat V_i$ and $\widehat V_i$
is tangent to $\widehat f^1(V)$ along $\Sigma_i$. Here $S_i$ is a wrinkle
of the wrinkled embedding $\widehat f^1$. The foliation $\F$ is transversal to
$\widehat V_i$ near $\Sigma_i$. Each wrinkle originates from
embryo (i.e. from a point) and hence one can rotate the tangent $(n-1)$-planes
to $\Sigma_i\subset \widehat V_i$ in $\widehat V_i$ to a position transversal
to the foliation $\F\cap \widehat V$. Such a rotation can be approximated
near $\Sigma'_i=\widehat f^1(S'_i)$  by an isotopy
of $\Sigma_i$ in $\widehat V_i$. Then, by our {\it inductional hypothesis}
there exists an  isotopy of $\Sigma_i$ in $\widehat V_i$, fixed near $\Sigma_i'$
and
such that the singularities of the final embedding
with respect to $\F\cap \widehat V_i$ are generalized $\,(n-1,q)$-wrinkles.
The resulting isotopies $g^t_i$, $t\in [1,2]$ (in $\widehat V_i$\,, for all $i$)
can be extended to an isotopy $\widehat f^t$, $t\in [1,2]$,
of the wrinkled embedding $\widehat f^1$.
Finally, we regularize the family $\widehat f^t$, $t\in [0,2]$,
in such a way that for $\widehat f^2$ our regularization is the 
$\F$-regularization as in Lemma \ref{thm:F-regularization-3},
and get the required family $f^\tau$, $\tau\in [0,1], \tau=2t$,
of embeddings such that  $f^1$ is a generalized wrinkled map with
respect to $\F$.
\qed

\btl\,
{\bf Remark.}
Using in the above proof  Theorem  \ref{thm:int-tang-global-folds}
instead of \ref{thm:int-tang-global} we get    a version of Theorem 
\ref{thm:embedding-foliation}
about embeddings with the double fold type tangency singularities
to a foliation.
\btr

\subsection{Embeddings into distributions}\label{sec:embed-to-dis}

In this section we sketch a generalization of Theorem 
\ref{thm:embedding-foliation}
to the case of distributions. 

\mn
Let $\xi$ be a distribution on a manifold $W^m$, $\codim\,\xi=q$.
An embedding $f:V^n \to (W^m,\xi)$, \,$n \geqslant \codim\,\xi$,
is called {\it transversal to} $\xi$, if the
{\it reduced differential}
$$d^{\,\xi} f:TV \mathop{\longrightarrow}\limits^{df} TW
\mathop{\longrightarrow}\limits^{\pi^{\xi}} TW/\xi$$
is surjective. The {\it non-transversality} defines a variety 
$\Sigma_\xi$ of the $1$-jet space $J^1(V,W)$. For a general 
non-integrable $\xi$ one cannot define fold and cusp type tangency
through normal forms. However, the original Whitney-Thom definition
is applicable in this situation as well.

\n 
We say that $f$ has  at a point $p\in V$  a tangency to $\xi$ of 
{\it fold} type if
\begin{itemize}
\item $\corank\, d^{\,\xi}_p =1$.
\item $J^1(f):V\to J^1(V,W)$ is transverse to $\Sigma$;
\item $d^{\,\xi}_pf|_{T_p\Sigma (f)}:T_p\Sigma(f)\to TW_{f(p)}/\xi$
is injective.
\end{itemize}
For a codimension $1$ distribution $\xi$ defined by a Pfaffian equation
$\alpha=0$, the fold tangency point of an embedding $f:V\to W$ to $\xi$
are isolated non-degenerate zeroes of the induced $1$-form $f^*\alpha$ on $V$.

\n
Similarly, we define a tangency of {\it cusp} type by the same first
two conditions and, in addition by requiring that
$d^{\,\xi}_pf|_{T_p\Sigma (f)}:T_p\Sigma(f)\to TW_{f(p)}/\xi$
is not injective and  the $2$-jet section $J^2(f):V\to J^2(V,W)$
to be transversal to the singularity in $J^2(V,W)$ defined by the
non-injectivity condition.

\n
Combining the following observations, one can associate  with a co-oriented
fold an {\it index}, as in the case of a foliation. 
\bp
{\bf (Index for fold-type tangency to distribution)}
\bi
\item [\bf A.] Let $J^2(\xi)\to W^m$ be the bundle of $\,2$-jets
of codimension $q$ submanifolds of $W$ tangent to $\xi$. There exists a section
$W\to J^2(\xi)$ which can be uniquely characterized by the following property:
for every point $p\in W$ there exists an embedding 
$\varphi_p:\bbR^{m-q}\to W$ such that $d_p\varphi_p(\bbR^{m-q})=\xi_p$,
and for all $t\in \bbR,\,x\in\bbR^{m-q}$ 
we have $\frac{d\varphi_p(tx)}{dt}\in\xi$.
\item [\bf B.] Let $(W,\xi)$ be as above, and $\Sigma$ a submanifold of
dimension $<\codim\,\xi$ such that $\pi^\xi:T\Sigma\to TW/\xi$ is injective.
Then, given any  section $\sigma:\Sigma\to J^2(\xi)$ there exists
codimension $q$ foliation $\F$ on $\Op\Sigma$ whose $2$-jet along $\Sigma$
is equal to  $\sigma$.
\item [\bf C.] Suppose that an embedding $f:V\to W$ has a fold type tangency of
index $k$ to a foliation $\F$ along a submanifold $\,\Sigma\subset W$.
Let $\wt \F$ be another foliation of the same codimension and which
have the same $2$-jet as $\F$ along $\Sigma$. Then  $f:V\to W$ has at
$\Sigma$ a fold type tangency to a foliation $\wt \F$  of the same index $k$. 
\ei
\ep
This allows us define generalized $(n,q)$-wrinkles and embryo type
singularities of tangency of an embedding $f:V^n\to (W,\xi)$ 
to a distribution $\xi$.  

\n
With these definitions Theorem \ref{thm:embedding-foliation}
can be generalized without any changes to the case
of an arbitrary distribution $\xi$ instead of a foliation.
The inductional proof for the foliation case presented in
Section \ref{sec:embed-to-fol} works in this more general with
a few additional remarks. Namely, one note that the base of the induction,
i.e. the case when $\dim\,\xi=1$ is the same in this case, as $1$-dimensional
distributions are integrable. When applying in the inductional step
the standard regularization from \ref{sec:wr-embed-formal}\,B,
one need to choose the regularizing function $\wt u(y,z)$ 
sufficiently $C^1$-close to $u(y,z)$.
With this modification the proof goes through as is.

%\bp
%{\bf (Corollary)}
%Let $\xi=\{\alpha=0\}$ be a codimension $1$ distribution, 
%$B$ an $(n+1)$-dimensional ball,  
%and $f_p:B \to W, p\in B$,  a family of embeddings $B\to W$.
%Then there exists a $C^0$-small isotopy $f^t_p$, $t\in[0,1]$,
%of the family $f^0_p=f_p$, $p\in B$, such that the family of 
%induced $1$-forms $\alpha_p=\left( f^1_p|_{\p B}\right)$
%have only non-degenerate, or birth-death zeroes. If for 
% over a closed subse $B_0\subset B$ the zeroes of the forms 
%$\alpha_p$ are non-degenerate, that the isotopy can be chosen fixed over $B_0$.
%\ep

\section{Appendix: Wrinkling}\label{appendix}

We recall here, for a convenience of the reader,
some definitions from \cite{[EM97]} and \cite{[EM98]}
and introduce the notions of {\it generalized wrinkles} and
{\it generalized wrinkled maps}. We also formulate here 
some {\it results} from \cite{[EM97]} and \cite{[EM98]}
though we do not use these theorems in the paper. 

\subsection{Folds and cusps}\label{sec:fold-cusp}

Let $V$ and $W$ be smooth manifolds of dimensions $n$ and $q$, respectively,
and $n\geq q$.  For a smooth map $f:V \rightarrow W$ we will denote by
$\Sigma(f)$ the set of its singular points, i.e.
$$
\Sigma(f) = \left\{ p \in V, \; \hbox { rank } d_pf < q \right\}\;.
$$
A point $p \in \Sigma(f)$ is called a {\it fold} type singularity or
a {\it fold} of index $s$ if near the point $p$ the map $f$ is equivalent to
the map
$$
\bbR ^{q-1} \times \bbR ^{n-q+1} \rightarrow \bbR ^{q-1}
\times \bbR ^{1}
$$
given by the formula
$$
(y,x) \mapsto \left( y,\,\, -\sum^s_1 x^2_i + \sum_{s+1}^{n-q+1} x^2_j
\right) \eqno(5)
$$
where $x=(x_1, \dots, x_{n-q+1}) \in \bbR ^{n-q+1}$ and $y=(y_1,...,y_{q-1})
\in \bbR ^{q-1}$.
For $W=\bbR ^1$ this is just a nondegenerate index $s$
critical point of the function $f:V \rightarrow \bbR ^1$.

\mn
Let $q>1$. A point $p \in \Sigma(f)$ is called a
{\it cusp} type singularity or a {\it
cusp} of index $s + \frac{1 }{ 2}$\; if near the point $p$ the map $f$ is
equivalent to the map
$$
\bbR ^{q-1}\times \bbR ^1\times \bbR ^{n-q}
\rightarrow \bbR^{q -1} \times \bbR ^1
$$
given by the formula
$$
(y,z,x) \mapsto \left( y,z^3 + 3y_1z - \sum^s_1 x_i^2 + \sum^{n-
q}_{s+1} x^2_j \right) \eqno (6)
$$
where $x=(x_1, \ldots, x_{n-q}) \in \bbR ^{n-q},\,\, z \in \bbR ^1,\,\,
y = (y_1, \ldots, y_{q-1}) \in \bbR ^{q-1}$.

\mn
For $q\geq 1$ a point $p \in \Sigma(f)$ is called an {\it embryo type}
singularity or an {\it embryo} of index $s + \frac{1}{ 2}$\,\, if $f$
is equivalent near $p$ to the map
$$
\bbR^{q-1}\times\bbR^1\times\bbR^{n-q}\rightarrow
\bbR^{q-1}\times\bbR ^1
$$
given by the formula
$$
(y,z,x) \mapsto \left( y, z^3 + 3|y|^2z - \sum^s_1x^2_i + \sum_{s+1}^{n-q}
x^2_j\right) \eqno (7)
$$
where $x \in \bbR ^{n-q}$,\,\, $y \in \bbR ^{q-1}$,\,\, $z \in \bbR ^1,\,\,
|y|^2=\sum\limits_1^{q-1}y_i^2$.

\mn
Notice that folds and cusps are stable singularities for individual
maps, while embryos are stable singularities only for $1$-parametric
families of mappings.  For a generic
perturbation of an individual map
embryos either disappear or give birth to wrinkles which we consider
in the next section.

\subsection{Wrinkles and wrinkled mappings}
\label{sec:wrinkles}

Consider the map
$$ w(n,q,s): \bbR^{q-1}  \times \bbR^1\times \bbR^{n-q}
\rightarrow \bbR^{q-1} \times \bbR^1 $$
given by the formula
$$ (y,z,x) \mapsto
\left( y, z^3 + 3(|y|^2-1)z - \sum^s_1 x^2_i + \sum_{s+1}^{n-q}
x_j^2 \right), \eqno (8) $$
where $y \in \bbR ^{q-1},\,\, z \in \bbR^1 ,\,\, x\in \bbR ^{n-q}$ and
$|y|^2 = \displaystyle{\sum_1^{q-1}y^2_{i}}$.

\mn
Notice that the singularity
$\Sigma (w(n,q,s))$ is the $(q-1)$-dimensional sphere
$$S^{q-1}=S^{q-1}\times 0\subset\bbR^q\times\bbR ^{n-q}.$$
Its equator
$\{ |y|=1, z=0, x=0\} \subset \Sigma (w(n,q,s))$
consists of cusp
points of index $s + \frac{1 }{ 2}$\,.  The upper hemisphere
$\Sigma (w) \cap\{z>0\}$ consists of folds of index $s$ and
the lower one $\Sigma (w) \cap \{z< 0\}$ consists
of folds of index $s+1$.
Also it is useful to notice that
the restrictions of the map $w(n,q,s)$ to subspaces $y_1=t$,
viewed as maps $\bbR ^{n-1} \to \bbR ^{q-1}$, are non-singular maps
for $|t|>1$, equivalent to $w(n-1,q-1,s)$ for $|t|<1$ and to embryos for
$t=\pm 1$.

\mn
Although the differential
$dw(n,q,s):T(\bbR^n) \rightarrow T(\bbR^q)$
degenerates at points of $\Sigma(w)$, it can be canonically {\it regularized}
over $\Op_{\bbR^n}D^q$, an open neighborhood of the disk
$D^q=D^q\times 0\subset \bbR^q\times\bbR^{n-q}$.
Namely, we can change the element $3(z^2 + |y|^2 -1)$ in the Jacobi matrix of
$w(n,q,s)$ by a function $\gamma$ which
coincides with $3(z^2+|y|^2-1)$ on $\bbR^n\setminus\Op_{\bbR^n}D^q$ and does
not vanish along the $q$-dimensional subspace
$\{x=0\}=\bbR^{q}\times {\bf 0}\subset \bbR^n\,.$
The new bundle map
${\cal R}(dw):T (\bbR ^n) \rightarrow T(\bbR ^q)$
provides a homotopically canonical
extension of the map
$dw:T(\bbR ^n\setminus \Op_{\bbR^n}D^q) \rightarrow T(\bbR ^q)$ to
an epimorphism (fiberwise surjective bundle map)
$T(\bbR ^n) \rightarrow T(\bbR ^q)$.  We call
${\cal R}(dw)$ the {\it regularized differential} of the map $w(n,q,s)$.

\mn
A smooth map $f:V^n\to W^q$, $n\geqslant q$, is called {\it wrinkled}, if
any connected component $S_i$ of the singularity
$\Sigma(f)$ is diffeomorphic
to the standard $(q-1)$-dimensional sphere $S^{q-1}$
and bounds in $V$ a $q$-dimensional disk $D_i$, such that
the map $f|_{\Op D_i}$ is equivalent to the map $w(n,q,s)|_{\Op_{\bbR^n}D^q}$.
The spheres $S_i$ and its images $f(S_i)$ are called {\it wrinkles} of the
wrinkled mapping $f$. The differential $df:T(V) \rightarrow T(W)$ of
the wrinkled map $f$ can be regularized (near each wrinkle and hence globally)
to obtain an epimorphism $\R (df):T(V) \rightarrow T(W)$.

\subsection{Fibered wrinkles and fibered wrinkled mappings}
\label{sec:fibrinkles}

For any integer $k>0$ the map $w(k+n,q,s)$
can be considered as a fibered map over
$\bbR^k\times {\bf 0}\subset \bbR^{k+n}$.
We shall refer to this fibered map as $w_k(k+n,q,s)$.
The regularized differential ${\cal R}(dw_k(k+n,q,s))$ is a
fibered (over $\bbR^k$) epimorphism
$$
\xymatrix
{\bbR^k \times T(\bbR^{q-1}\times\bbR^1\times \bbR^{n-q})
\ar[rrr]^{\,\,\,\,\,\,\,\,{\cal R}(dw_k(k+n,q,s))}&&&\bbR^k
\times T(\bbR^{q-1} \times \bbR^1)}
$$
A fibered (over $B$) map $f:V^{k+n}\to W^{k+q}$, $n\geqslant q$, is called
{\it fibered  wrinkled}, if any connected component $S_i$
of the singularity $\Sigma(f)$ is diffeomorphic
to the standard $(k+q-1)$-dimensional sphere $S^{k+q-1}$
and bounds in $V$ a $(k+q)$-dimensional disk $D_i$, such that
the fibered map $f|_{\Op D_i}$ is equivalent to the
fibered map $w_k(k+n,q,s)|_{\Op_{\bbR^{k+n}}D^q}$.
The spheres $S_i$ and its images $f(S_i)$ are called {\it fibered wrinkles}
of the fibered wrinkled mapping $f$. The restrictions of
a fibered wrinkled map to a fiber may have, in addition to wrinkles,
emdryos singularities. For a fibered wrinkled map $f:V\to W$
one can define its regularized differential
which is a fibered (over $B$) epimorphism
$\R(d_Bf):T_BV\to T_BW$.

\subsection{Main theorems about wrinkled mappings}
\label{sec:main-wrinkling}

The following  Theorem \ref{thm:A}
and its parametric version \ref{thm:B}
are the main results of our paper \cite{[EM97]}:
\bp
[Wrinkled mappings]\label{thm:A}
Let $F: T(V) \rightarrow T(W)$ be an
epimorphism which covers a map $f: V\rightarrow W$.
Suppose that $f$ is a submersion on a neighborhood of a
closed subset $K \subset V$, and $F$
coincides with $df$ over that neighborhood.
Then there exists a wrinkled map $g : V \rightarrow W$ which coincides with
$f$ near $K$ and such that ${\cal R} (dg)$ and $F$ are homotopic rel. $T(V)|_K$.
Moreover, the map $g$ can be chosen arbitrarily $C^0$-close to $f$,
and his wrinkles can be made arbitrarily small.
\ep

\bp
[Fibered wrinkled mappings]\label{thm:B}
Let $f: V \rightarrow W$ be a fibered over $B$ map covered by a fibered
epimorphism $F: T_B(V) \rightarrow T_B(W)$.  Suppose that $F$ coincides with
$df$ near a closed subset $K \subset V$ (in particular, $f$ is a fibered
submersion near $K$), then there exists a fibered wrinkled map
$g: V\rightarrow W$
which extends $f$ from a neighborhood of $K$, and such that the
fibered epimorphisms ${\cal R}(dg)$ and $F$ are homotopic
rel. $T_B(M)|_K$. Moreover, the map $g$ can be chosen arbitrarily $C^0$-close
to $f$, and his wrinkles can be made arbitrarily small.
\ep

\subsection{Wrinkled mappings into foliations}
\label{sec:mappings into foliations}

Here we formulate a slightly strengthened version of  Theorems \ref{thm:A}
and \ref{thm:B}. Let us start with some definitions.

\mn
Let $\F$ be a foliation  on  a  manifold $W$, $\codim\,\F=q$.
A map $f:V \to W$ is
called {\it transversal to} $\F$, if the reduced differential
$$
TV \mathop{\longrightarrow}\limits^{df} TW
 \mathop{\longrightarrow}\limits^{\pi_{\F}} \nu(\F)
$$
is an epimorphism. Here $\nu(\F)=TW/\tau(\F)$ is the normal bundle of
the foliation $\F$.

\mn
An open subset $U\subset W$ is called
{\it elementary} (with respect to $\F$), if $\F|_U$ is generated by
a submersion $p_{_U}:U\to \bbR^q$.
An open subset
$U\subset V$ is called {\it small} (with respect to $f$ and $\F$),
if $f(U)$  is contained in an elementary subset $U'$ of $W$.
A map $f:V\to W$ is called $\F^{\bot}$-{\it wrinkled}, or {\it wrinkled
with respect to $\F$}, if there exist
disjoint small subsets $U_1,...U_l\subset V$ such that
$f|_{V\setminus (U_1\cup...\cup U_l)}$ is transversal to $\F$
and for each $i=1,...,l$ the composition
$$
U_i \mathop{ \longrightarrow}\limits^{f|_{U_i}} U'_i
\mathop{ \longrightarrow}\limits^{p_{{U'_i}}} \bbR^q
$$
(where $U'_i\supset f(U_i)$ is an elementary subset of $W$),
is a wrinkled map.
In order to get the {\it regularized} reduced differential
$$
\R(\pi_\F\circ df):TM\to \nu(\F)
$$
of the $\F^{\bot}$-wrinkled map $f$, 
we regularize the differential of each  wrinkled map
$w_i=p_{_{U'_i}}\circ f|_{U_i}$ as in Sections \ref{sec:wrinkles}
and then set
$$\R(\pi_\F\circ df|_{U_i})=[dp_{U'_i}|_{\nu(\F)}]^{-1}\circ\R(dw_i).$$

\mn
Similarly to Section \ref{sec:fibrinkles} we can define
a {\it fibered} $\F^\bot$-wrinkled map $f:V\to Q$,
where the foliation $\F$ on $W$ is fibered over
the same base $B$. Finally we define, in a usual way,
the {\it regularization}
$$
\R(\pi_\F\circ d_Bf):T_BV\to \nu_B(\F)
$$
of the fibered reduced differential
$$
T_BV \mathop{\longrightarrow}\limits^{d_Bf} TW
\mathop{\longrightarrow}\limits^{\pi_{\F}} \nu_B(\F).
$$

\bp
\label{thm:Af}
{\bf (Wrinkled mappings of manifolds into foliations)}
Let $\F$ be a foliation on a manifold $W$ and let
$F:TV\to \nu(\F)$ be an epimorphism which covers a map
$f:V\to W$. Suppose that $f$ is transversal to $\F$ in
a neighborhood of closed subset $K$ of $V$, and $F$ coincides
with the reduced differential $\pi_\F\circ df$ over that neighborhood.
Then there exists a  $\F^\bot$-wrinkled map $g:V\to W$ such that
$g$ coincides with $f$ near $K$, and
 $\R(\pi_\F\circ dg)$ is homotopic to $F$ relative to  $TV|_K$.
\ep

\n
{\bf Proof.}
Take a triangulation of the manifold $V$
by {\it small} simplices. First we use Gromov-Phillips'
theorem (see \cite{[Ph67]}, \cite{[Gr86]} or \cite{[EM02]}) to approximate
$f$  near the $(n-1)$-skeleton of the triangulation by a map transversal
to $\F$. Then, using Theorem \ref{thm:A},
for a neighborhood $U_i$ of every $n$-simplex $\sigma_i$
and an elementary set $U'_i\supset f(U_i)$ we can approximate the map
$p_{U'_i}\circ f|_{U_i}$ by a  wrinkled map.
This approximation can be realized by a deformation of the map $f$,
keeping it fixed on a closed subset of $U_i$, where $f$
was already previously defined. This process produces
the desired  $\F^\bot$-wrinkled map.
$\square$

\mn
Similarly, Theorem \ref{thm:B} can be generalized to
the following fibered version of Theorem \ref{thm:Af}.

\bp
\label{thm:Bf}
{\bf (Fibered wrinkled mappings of manifolds into foliations)}
Let $f: V \rightarrow W$ be a fibered over $B$ map, $\F$ be a fibered over $B$
foliation on $W$ and let $F: T_B(V) \rightarrow \nu_B(\F)$ be a fibered
epimorphism which covers $f$. Suppose that $f$ is
fiberwise transversal to $\F$  near a closed subset
$K \subset V$, and  $F$ coincides with fibered reduced differential
$$
T_BV \mathop{\longrightarrow}\limits^{d_Bf} TW
 \mathop{\longrightarrow}\limits^{\pi_{\F}} \nu_B(\F).
$$
near $K$.
Then there exists a fibered  $\F^\bot$-wrinkled map
$g: V\rightarrow W$ which extends $f$ from a neighborhood of $K$, and such that
the fibered epimorphisms ${\cal R}(\pi_\F\circ d_Bg)$ and $F$ are
homotopic rel. $T_B(V)|_K$.
\ep

\subsection{Generalized wrinkled mappings}
\label{sec:generalized wrinkles}

A smooth map $f:V^n\to W^q$, $n\geqslant q$, is called
{\it generalized wrinkled}, if any connected components $S_i$
of the singularity $\Sigma=\Sigma_f$ is diffeomorphic
to the standard $(q-1)$-dimensional sphere $S^{q-1}$, which bounds in $V$
a $q$-dimensional disk $D_i$,  and for each such sphere
the map $f|_{\Op S_i}$ is equivalent to the map
$w(n,q,s)|_{\Op_{\bbR^n}{S^{q-1}}}$.
The spheres $S_i$ and its images $f(S_i)$ are called {\it generalized wrinkles}
of the generalized wrinkled mapping $f$.

\mn
A fibered map $f:V^{k+n}\to W^{k+q}$, $n\geqslant q$,
is called {\it generalized fibered wrinkled map}, if
any connected components $S_i$ of the singularity
$\Sigma=\Sigma_f$ is diffeomorphic
to the standard $(k+q-1)$-dimensional sphere $S^{k+q-1}$, which bounds in $V$
a $(k+q)$-dimensional disk $D_i$, and for each such sphere
the fibered map $f|_{\Op S_i}$ is equivalent to
the fibered map $w_k(k+n,q,s)|_{\Op_{\bbR^{k+n}}{S^{q-1}}}$.
The spheres $S_i$ and its images $f(S_i)$ are called {\it generalized
fibered wrinkles} of the generalized wrinkled mapping $f$.

\mn
Let $\F$ be a foliation  on  a  manifold $W$, $\codim\,\F=q$.
A map $f:V\to W$ is called {\it generalized wrinkled
with respect to $\F$}, 
if there exist disjoint small subsets $U_1,...U_l\subset V$ such that
$f|_{V\setminus (U_1\cup...\cup U_l)}$ is transversal to $\F$
and for each $i=1,...,l$ the composition
$$
U_i \mathop{ \longrightarrow}\limits^{f|_{U_i}} U'_i
\mathop{ \longrightarrow}\limits^{p_{{U'_i}}} \bbR^q
$$
(where $U'_i\supset f(U_i)$ is an elementary subset of $W$),
is a generalized wrinkled map.

%\clearpage

\end{document}